\documentclass{nyj}

\newcommand{\floor}[1]{\lfloor #1 \rfloor }

\def\R{{\hbox{\bf R}}}
\def\P{{\hbox{\bf P}}}
\def\Z{{\hbox{\bf Z}}}
\def\C{{\hbox{\bf C}}}
\def\Q{{\hbox{\bf Q}}}

\def\B{{\hbox{\bf B}}}
\def\I{{\hbox{\bf I}}}

\def\A{{\hbox{\bf A}}}

\def\Ca{{K}}
\def\Cc{{C_1}}
\def\Cb{{C_2}}
\def\Cd{{C_3}}
\def\Ce{{K}}
\def\Cf{{C_5}}
\def\Cg{{C_6}}
\def\Ch{{C_7}}
\def\Ci{{C_8}}
\def\Cj{{C_9}}

\def\Cl{{C_{11}}} 

\font \roman = cmr10 at 10 true pt
\def\x{{\roman x}}
\def\y{{\roman y}}

\def\be#1{\begin{equation}\label{#1}}
\def\bas{\begin{align*}}
\def\eas{\end{align*}}
\def\bi{\begin{itemize}}
\def\ei{\end{itemize}}

\def\dim{{\hbox{\roman dim}}}

\def\dist{{\hbox{\roman dist}}}

\def\eps{\varepsilon}
\def \endprf{\hfill  {\vrule height6pt width6pt depth0pt}\medskip}
\def\emph#1{{\it #1}}
\def\textbf#1{{\bf #1}}
\theoremstyle{plain}
  \newtheorem{theorem}[subsection]{Theorem}

  \newtheorem{proposition}[subsection]{Proposition}
  
  \newtheorem{lemma}[subsection]{Lemma}
  \newtheorem{corollary}[subsection]{Corollary}
  \newtheorem{Distance conjecture}[subsection]{Distance Conjecture}
  \newtheorem{Bilinear distance conjecture}[subsection]{Bilinear Distance Conjecture}
  \newtheorem{Furstenburg problem}[subsection]{Furstenburg problem}
  \newtheorem{Discretized Furstenburg conjecture}[subsection]{Discretized Furstenburg 
  Conjecture}
  \newtheorem{Ring problem}[subsection]{Ring Problem}
  \newtheorem{Ring conjecture}[subsection]{Ring Conjecture}
  \newtheorem{Main theorem}[subsection]{Main Theorem}
  \newtheorem{Cauchy-Schwarz}[subsection]{Cauchy-Schwarz}
  \newtheorem{Refinement}[subsection]{Refinement}
  \newtheorem{Separation}[subsection]{Separation}
  \newtheorem{Perfection}[subsection]{Perfection}
  \newtheorem{Kakeya}[subsection]{Kakeya}

\theoremstyle{remark}

\theoremstyle{definition}
  \newtheorem{definition}[subsection]{Definition}

\include{psfig}

\begin{document}

\title[Falconer and Furstenburg]{Some connections between Falconer's distance set 
conjecture, and sets of Furstenburg type }

\author{Nets Hawk Katz}
\address{Department of Mathematics, University of Illinois at Chicago, Chicago IL 60607-7045}
\email{nets@@math.uic.edu}

\author{Terence Tao}
\address{Department of Mathematics, UCLA, Los Angeles CA 90095-1555}
\email{tao@@math.ucla.edu}


\begin{abstract} 
In this paper we investigate three unsolved conjectures in
geometric combinatorics, namely Falconer's distance set conjecture,
the dimension of Furstenburg sets, and Erd\"os's ring conjecture.
We formulate natural $\delta$-discretized versions of these conjectures
and show that in a certain sense that these discretized versions are
equivalent.
\end{abstract}

\maketitle
\tableofcontents

\section{Introduction}

In this paper we study Falconer's distance problem, the dimension of sets
of Furstenburg type, and Erd\"os's ring problem.
Although we have no direct progress on any of these problems, we are
able to reduce the geometric problems to $\delta$-discretized variants
and show that these variants are all equivalent.

In order to state the main results we first must develop a certain amount
of notation.

\subsection{Notation}

$0 < \eps \ll 1$, $0 < \delta \ll 1$ are small parameters.  We use $A \lessapprox B$ to 
denote the estimate $A \leq C_\eps \delta^{-C\eps} B$ for some constants $C_\eps$, $C$, 
and $A \approx B$ to denote $A \lessapprox B \lessapprox A$.  

We use $\B(x,r) = \B^n(x,r)$ to denote the open ball of radius $r$ centered at $x$ in 
$\R^n$, and $\A = \A^n$ to denote any annulus in $\R^n$ of the form 
$\A := \{ x: |x| \approx 1 \}$.

If $A$ is a finite set, we use $\# A$ to denote the cardinality of $A$.
For finite sets $A$, $B$, we say that $A$ is a refinement of $B$ if $A \subset B$ and 
$\# A \approx \# B$.

If $E$ is contained in a subspace of $\R^n$ and has positive measure in that subspace, 
we use $|E|$ for the induced Lebesgue measure of $E$.  The subspace will always be clear 
from context.

For sets $E, F$ of finite measure, we say that $E$ is a refinement of $F$ if $E \subset F$ 
and $|E| \approx |F|$.
We say that $E$ is \emph{$\delta$-discretized} if $E$ is the union of balls of radius 
$\approx \delta$.

\begin{definition}  For any $0 < \alpha \leq n$, we say that a set $E$ is a 
$(\delta,\alpha)_n$-set if it is contained in a ball $\B^n(0,C)$, is 
$\delta$-discretized and one has
\be{dim}
|E \cap \B(x,r)| \lessapprox \delta^n (r/\delta)^\alpha
\end{equation}
for all $\delta \leq r \leq 1$ and $x \in \R$. 
\end{definition}

Roughly speaking, a $(\delta,\alpha)_n$-set behaves like the $\delta$-neighbourhood of an 
$\alpha$-dimensional set in $\R^n$.  The condition \eqref{dim} is necessary to ensure that 
$E$ does not concentrate in a small ball, which would lead to some trivial counterexamples 
to the conjectures in this paper.  (cf. the ``two ends'' condition in \cite{wolff:kakeya},
\cite{wolff:xray}).

If $X, Y$ are subsets of $\R^n$, we use $X+Y$ to denote the set $X+Y := \{ x+y: x \in X, 
y \in Y\}$.  Similarly define $X - Y$, and (when $n=1$) $X \cdot Y$, $X/Y$, $X^2$, 
$\sqrt{X}$, etc.  Note that $X^2 \subsetneq X \cdot X$ in general.  Note that $X \times Y$
 denotes the Cartesian product $X \times Y := \{ (x,y): x \in X, y \in Y \}$ as opposed to 
the pointwise product $X \cdot Y := \{ xy: x \in X, y \in Y\}$.  Unfortunately there is a 
conflict of notation between $X^2 := \{ x^2: x \in X \}$ and $X^2 := \{ (x,y): x,y \in X \}$; 
to separate these two we shall occasionally write the latter as $X^{\oplus 2}$.

If a rectangle $R$ has sides of length $a, b$ for some $a > b$, we call the \emph{direction} 
of $R$ the direction $\omega \in S^1$ that the sides of length $a$ are oriented on.  This is 
only defined up to sign $\pm$.

\subsection{The Falconer distance problem}

For any compact subset $K$ of the plane $\R^2$, define the \emph{distance set} $\dist(K) 
\subset \R$ of $K$ by
$$ \dist(K) := |K-K| = \{ |x-y|: x,y \in K \}.$$
In \cite{falconer:distance} Falconer conjectured that if $\dim(K) \geq 1$, then
$\dim(\dist(K)) = 1$, where $\dim(K)$ denotes the Hausdorff dimension of $K$. As progress 
towards this conjecture, it was shown in \cite{falconer:distance} that $\dim(\dist(K)) = 1$ 
obtained whenever $\dim(K) \geq 3/2$.  This was improved to $\dim(K) \geq 13/9$ by Bourgain 
\cite{borg:distance} and then to $\dim(K) \geq 4/3$ by Wolff \cite{wolff:distance}.  These 
arguments are based around estimates for $L^2$ circular means of Fourier transforms of 
Frostman measures.  However, it is unlikely that a purely Fourier-analytic approach will be 
able to improve upon the 4/3 exponent; for a discussion, see \cite{wolff:distance}.

Now suppose that one only assumes that $\dim(K) \geq 1$.  An argument of Mattila 
\cite{mattila:distance} shows that $\dim(\dist(K)) \geq \frac{1}{2}$.  One may ask whether 
there is any improvement to this result, in the following sense:

\begin{Distance conjecture}\label{dist-conj} There exists an absolute constant $c_0 > 0$ such that 
$\dim(\dist(K)) \geq \frac{1}{2} + c_0$ whenever $K$ is compact and satisfies 
$\dim(K) \geq 1$.
\end{Distance conjecture}

This is of course weaker than Falconer's conjecture, but remains open.

One may hope to prove this conjecture by first showing a $\delta$-discretized analogue.  
As a naive first approximation, we may ask the informal question of whether (for 
$0 < \delta, \eps \ll 1$) the distance set of a $(\delta,1)_2$ set of measure 
$\approx \delta$ can be (mostly) contained in a $(\delta,1/2)_1$ set.

Unfortunately, this problem has an essentially negative answer, as the counterexample
\be{counter}
\{ (x_1,x_2): x_1 = k \sqrt{\delta} + O(\delta), x_2 = O(\sqrt{\delta}) \hbox{ for some } 
k \in \Z, k = O(\delta^{-1/2}) \}
\end{equation}
shows\footnote{This counterexample also appears in Fourier-based approaches to the distance
problem, see \cite{wolff:distance}}.  A substantial portion of the distance set of 
\eqref{counter} is contained in the $\delta$-neighbourhood of an arithmetic progression of 
spacing $\delta^{1/2}$, and this is a $(\delta,\frac{1}{2})_1$ set.

 \begin{figure}[htbp] \centering
 \ \psfig{figure=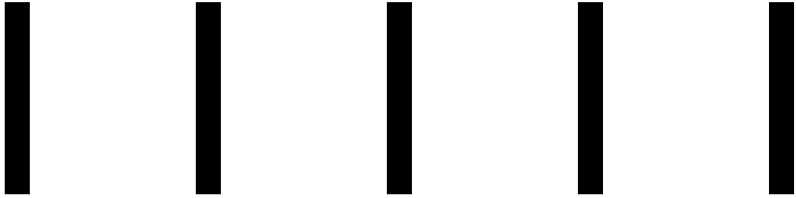}
 \caption{An example to remember. Few blurred distances but many blurred points.}
 \end{figure}
 
This obstruction to solving Conjecture \ref{dist-conj} can be eliminated by replacing the 
above informal problem with a ``bilinear'' variant in which an angular separation condition 
is assumed:\footnote{This idea is frequently used in related problems, see e.g. 
\cite{tvv:bilinear}, \cite{borg:cone}, \cite{wolff:cone}.  Other discretizations are 
certainly possible, providing of course that \eqref{counter} is neutralized.}

\begin{Bilinear distance conjecture}\label{less-naive}  Let $Q_0$, $Q_1$, $Q_2$ be three 
cubes in $\B(0,C)$ of 
radius $\approx 1$ satisfying the separation condition
\be{gen-sep}
 | (x_1-x_0) \wedge (x_2-x_0) | \approx 1 \hbox{ for all } x_0 \in Q_0, x_1 \in Q_1, 
x_2 \in Q_2.
\end{equation}
For each $j=0,1,2$, let $E_j$ be a $(\delta,1)_2$ subset of $Q_j$, and let $D$ be a 
$(\delta,1/2)_1$ subset of $\R$.  Then
\be{eps1}
|\{ (x_0,x_1,x_2) \in E_0 \times E_1 \times E_2: |x_0-x_1|, |x_0-x_2| \in D \}| \lessapprox 
\delta^{3-c_1}
\end{equation}
where $c_1 > 0$ is an absolute constant.
\end{Bilinear distance conjecture}

The estimate \eqref{eps1} is trivially true when $c_1 = 0$.  Also, if it were not for 
condition \eqref{gen-sep} one could easily disprove \eqref{eps1} for any $c_1 > 0$ by 
modifying \eqref{counter}.  Conjecture \ref{less-naive} is also heuristically plausible from 
analogy with results on the discrete distance problem such as Chung, Szemer\'edi and 
Trotter \cite{cst:erdos}.  We remark that the arguments in that paper require the 
construction of three 
cubes satisfying \eqref{gen-sep}, and involve the Szemer\'edi-Trotter theorem (which may be 
considered as a result concerning the discrete analogue of the Furstenburg problem).

 \begin{figure}[htbp] \centering
 \ \psfig{figure=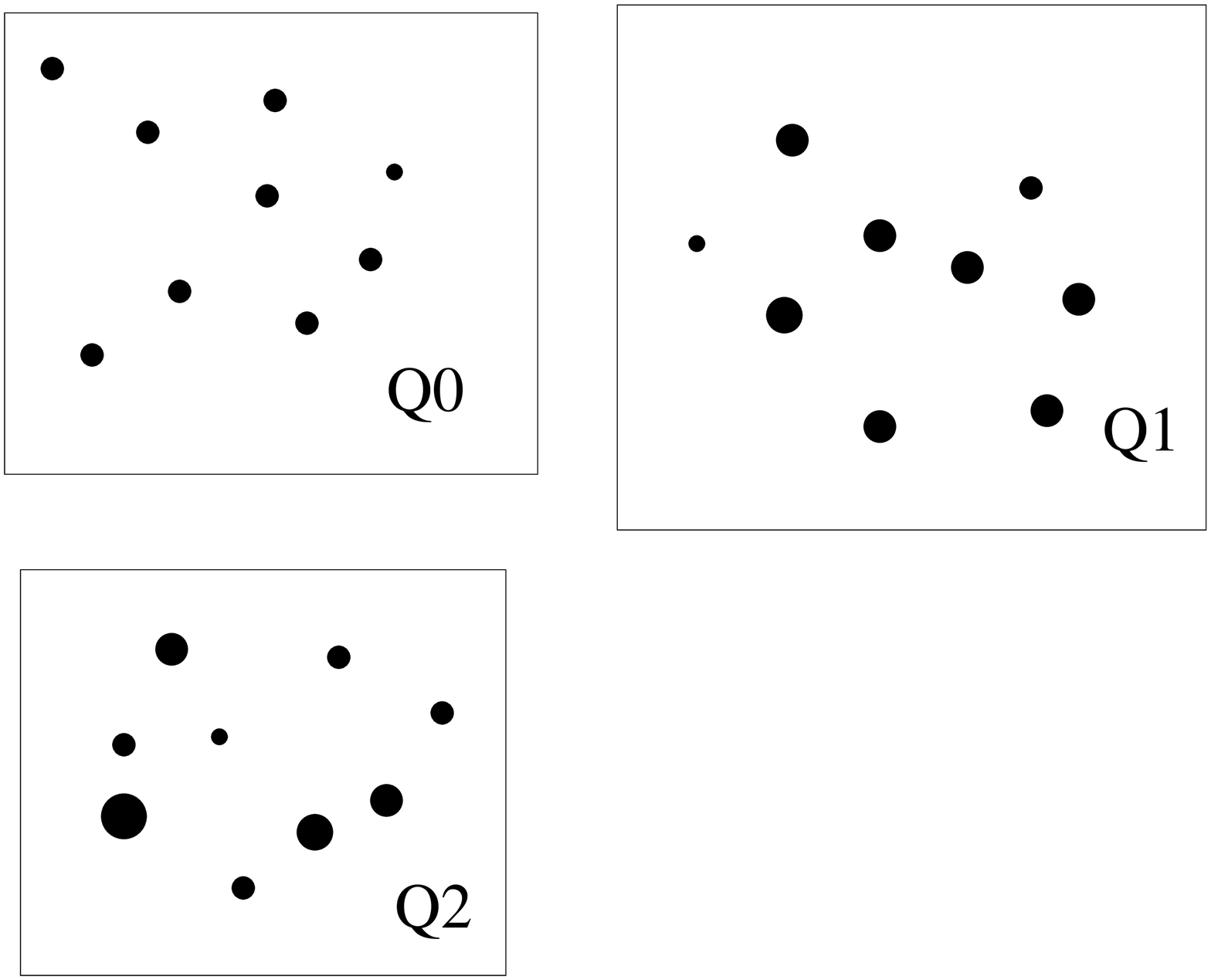,width=5in,height=5in}
 \caption{In the bilinear distance conjecture, the points are split into three camps.}
  \end{figure}

In Section \ref{falc-reduc-sec} we prove

\begin{theorem}\label{falc-reduc} A positive answer to the bilinear
distance conjecture \ref{less-naive} implies a 
positive answer to the distance conjecture \ref{dist-conj}.
\end{theorem}

Although this implication looks plausible from discretization heuristics, there are 
technical difficulties due to the presence of the counter-example \eqref{counter}, and also 
by the fact that several scales may be in play when studying the Hausdorff dimension of a 
set. 

\subsection{Dimension of sets of Furstenburg type}

We now turn to a problem arising from the work of Furstenburg, as formulated in work of 
Wolff \cite{wolff:survey}, \cite{wolff:distance}.  

\begin{definition}  Let $0 < \beta \leq 1$.  We define a $\beta$-set to be a compact set 
$K \subset \R^2$ such that for every direction $\omega \in S^1$ there exists a line segment 
$l_\omega$ with direction $\omega$ which intersects $K$ in a set with Hausdorff dimension 
at least $\beta$.  We let $\gamma(\beta)$ be the infimum of the Hausdorff dimensions of 
$\beta$-sets.
\end{definition}

In \cite{wolff:survey} the problem of determining $\gamma(\beta)$ is formulated.  At present 
the best bounds known are
$$ \max(\beta + \frac{1}{2}, 2\beta) \leq \gamma(\beta) \leq \frac{3}{2}\beta + 
\frac{1}{2},$$
see \cite{wolff:survey}.  This problem is clearly connected with the Kakeya problem 
(which is essentially concerned with the higher-dimensional analogue of $\gamma(1)$).  
Connections to the Falconer distance set problem have also been made, see 
\cite{wolff:distance}.

The most interesting value of $\beta$ appears to be $\beta = 1/2$.  In this case the two 
lower bounds on $\gamma(\beta)$ coincide to become $\gamma(\frac{1}{2}) \geq 1$.  We ask

\begin{Furstenburg problem}\label{furst-conj}  Is it true that $\gamma(\frac{1}{2}) 
\geq 1 + c_2$ for 
some absolute constant $c_2 > 0$?  In other words, is it true that $\frac{1}{2}$-sets must 
have Hausdorff dimension at least $1 + c_2$?
\end{Furstenburg problem}

One can $\delta$-discretize this problem as

\begin{Discretized Furstenburg conjecture}\label{furst-disc}  Let $0 < \delta \ll 1$, and 
let $\Omega$ be a 
$\delta$-separated set of directions, and for each $\omega \in \Omega$ let $R_\omega$ be a 
$(\delta,\frac{1}{2})_2$ set contained in a rectangle of dimensions $\approx 1 \times \delta$
 oriented in the direction $\omega$.  Let $E$ be a $(\delta,1)_2$ set.  Then
\be{furst-est} | \{ (x_0,x_1) \in E \times E: x_1, x_0 \in R_{\omega} \hbox{ for some } 
\omega \in \Omega \} |
\lessapprox \delta^{2+c_3}
\end{equation}
for some absolute constant $c_3 > 0$.
\end{Discretized Furstenburg conjecture}

As before, this conjecture is heuristically plausible from analogy with discrete incidence 
combinatorics, in particular the Szemer\'edi-Trotter theorem \cite{st:whatcanisay}.  Unlike 
the case with the 
distance problem, the set \eqref{counter} does not provide a serious threat, and so one does 
not need to go to a bilinear framework.  

The discretized Furstenburg conjecture \ref{furst-disc} is related to the
Furstenburg conjecture \ref{furst-conj} in much the same way 
that the Kakeya maximal function conjecture is related to the Kakeya set conjecture.  
In section \ref{furst-disc-sec} we show

\begin{theorem}\label{furst-equiv} A positive answer to the discretized
Furstenburg conjecture \ref{furst-disc} implies a 
positive answer to the Furstenburg problem \ref{furst-conj}.
\end{theorem}

\subsection{The Erd\"os ring problem}

We consider a problem of Erd\"os, namely

\begin{Ring problem} Does there exist a sub-ring $R$ of $\R$ which is a Borel set and has 
Hausdorff dimension strictly between 0 and 1?
\end{Ring problem}

This problem is connected to Falconer's distance problem; for instance, 
Falconer \cite{falconer:distance} used results on the distance problem to show that Borel 
sub-rings $R$ of $\R$ could not have Hausdorff dimension strictly between $1/2$ and 1.  
Essentially, the idea is to use the fact that $\dist(R \times R) \subseteq \sqrt{R}$.  

We concentrate on the specific problem of whether a sub-ring can have dimension exactly 1/2; 
it seems reasonable to conjecture that such rings do not exist.  
A positive answer to Conjecture \ref{dist-conj} would essentially imply this conjecture.

If $R$ is a ring of dimension $1/2$, then of course $R+R$ and $RR$ also have dimension 
$1/2$.  This leads us to the following $\delta$-discretization of the above conjecture.

\begin{Ring conjecture}\label{erdos}  Let $0 < \delta \ll 1$, and let $A \subset \A$ be a 
$(\delta,\frac{1}{2})_1$ set of measure $\approx \delta^{1/2}$.  Then at least one of 
$A+A$ and $AA$ has measure $\gtrapprox \delta^{\frac{1}{2}-c_4}$, where $c_4 > 0$ is an 
absolute constant.
\end{Ring conjecture}

The dimension condition \eqref{dim} is crucial, as the trivial counterexample 
$A := [1, 1 + \delta^{1/2}]$ demonstrates.  In principle the discretized ring
conjecture gives a negative answer to the Erd\"os ring problem, but we have
not been able to make this rigorous.

For the discrete version of this problem, when measure is replaced by cardinality, there is 
a result of Elekes \cite{Elekes:Gyorgy} that when $A$ has finite cardinality $\# A$, at least 
one of $A+A$ 
and $A A$ has cardinality $\gtrapprox \# A^{5/4}$ .  The proof of 
this result exploits the Szemer\'edi-Trotter theorem.  This is heuristic evidence for the 
ring conjecture \ref{erdos} if one accepts the 
(somewhat questionable) analogy between discrete models and 
$\delta$-discretized models.

It may appear that the ring hypothesis is being under-exploited when reducing to the ring 
conjecture \ref{erdos}, since one is only using the fact that $R+R$ and $RR$ are small.  
However, we shall see in Proposition \ref{cdcd} that control on $A+A$ and $AA$ actually 
implies quite good control on other arithmetic expressions such as $AA-AA$ or 
$(A-A)^2 + (A-A)^2$ (after passing to a refinement), so the ring hypothesis is not being 
wasted.

\subsection{The main result}

\centerline{\emph{One Ring to rule them all,}}
\centerline{\emph{One Ring to find them,}}
\centerline{\emph{One Ring to bring them all,}}
\centerline{\emph{and in the darkness bind them. \cite{tolkien}}}

As one can see from the previous discussion, there have been many partial connections drawn 
between the Falconer, Furstenburg, and Erd\"os problems.  The main result of this paper is 
to consolidate these connections into

\begin{Main theorem}\label{grand-equiv}  The conjectures \ref{less-naive}, \ref{furst-disc}, 
and \ref{erdos} are logically equivalent.
\end{Main theorem}

We shall prove this theorem in Sections \ref{arithmetic}-\ref{furst-falcon}.

In particular, in order to make progress on the Falconer and Furstenburg problems it 
suffices to prove the ring conjecture \ref{erdos}.  This appears to be the easiest of all 
the above problems to attack.
It seems likely that one needs to exploit some sort of ``curvature'' between addition and 
multiplication to prove this conjecture, although a naive Fourier-analytic pursuit of this 
idea seems to run into difficulties.  This may indicate that a 
combinatorial approach will be more fruitful than a Fourier approach.  The fact that $\R$ is 
a totally ordered field may also be relevant, since the analogue of Erd\"os's ring problem 
is false for non-ordered fields such as the complex numbers $\C$ or the finite field 
$F_{p^2}$.  (Unsurprisingly, the analogues of Falconer's distance problem and the 
conjectures for Furstenburg sets also fail for these fields, see e.g. 
\cite{wolff:survey}).

These problems are also related to the Kakeya problem in three dimensions, although the 
connection here is more tenuous.  A proof of Conjecture \ref{erdos} would probably lead 
(eventually!) to an alternate proof of the main result in \cite{katzlaba}, namely that 
Besicovitch sets\footnote{A Besicovitch set is a set which contains a unit line segment in 
every direction.} in $\R^3$ have Minkowski dimension strictly greater than 5/2, and would 
not rely as heavily
on the assumption that the line segments all point in different directions.  
Very informally, the point is that the arguments in \cite{katzlaba} can be pushed a bit 
further to conclude that a Besicovitch set of dimension exactly 5/2 must essentially be a 
``Heisenberg group'' over a ring of dimension 1/2.  We shall not pursue this connection in 
detail as it is somewhat lengthy and would not directly yield any new progress on the 
Kakeya problem.

In conclusion, these results indicate that the possibility of $1/2$-dimensional rings is a 
fundamental obstruction to further progress on the Falconer and Furstenburg problems, and 
may also be obstructing progress on the Kakeya conjecture and related problems (restriction, 
Bochner-Riesz, Stein's conjecture, local smoothing, etc.)  It also appears that 
substantially new techniques are needed to tackle this obstruction, possibly exploiting the 
ordering of the reals.

\section{Basic tools}

In this section $0 < \eps \ll 1$ is fixed, but $\delta$ is allowed to vary.  As in other 
sections, the implicit constants here are not allowed to depend on $\delta$.

To clarify many of the arguments in this paper, it may help to know that almost all 
estimates of the form $A \gtrapprox B$ which occur in this paper are sharp in the sense that 
the converse bound $A \lessapprox B$ is usually trivial to prove.  It is this sharpness 
which allows us to pass from one expression to another without losing very much in the 
estimates (if one does not mind the implicit constants in the $\lessapprox$ notation 
increasing very quickly!).

A typical application of this philosophy is

\begin{Cauchy-Schwarz}\label{cauchy}  Let $A$, $B$ be sets of finite measure, and let $\sim$ be a 
relation between elements of $A$ and elements of $B$.  If
$$ |\{ (a,b) \in A \times B: a \sim b \}| \geq \lambda |A| |B|$$
for some $0 < \lambda \leq 1$ then
$$ |\{ (a,b,b') \in A \times B \times B: a \sim b, a \sim b' \}| \geq \lambda^2 |A| |B|^2$$
\end{Cauchy-Schwarz}

\begin{proof}  We can rewrite the hypothesis as
$$ \int_A |\{ b \in B: a \sim b\}|\ da \geq \lambda |A||B|$$
and the conclusion as
$$ \int_A |\{ b \in B: a \sim b\}|^2\ da \geq \lambda^2 |A||B|^2.$$
The claim then follows from Cauchy-Schwarz.
\end{proof}

The next lemma deals with the issue of how to refine a $\delta$-discretized set to become a 
$(\delta,\alpha)_n$ set for suitable $\alpha$.

\begin{Refinement}\label{refine}  Let $0 < \delta \ll 1$ be a dyadic number, 
$0 < \alpha < n$, $K \gg 1$ be a constant, and let $E$ be a $\delta$-discretized set in 
$\B^n(0,C)$ such that $|E| \lessapprox \delta^{n-\alpha}$.  Then one can find a set 
$E_{\delta'}$ for all dyadic $\delta < \delta' \leq 1$ which can be covered by 
$\lessapprox \delta^{K\eps} {\delta'}^{-\alpha}$ balls of radius $\delta'$, and a set 
$(\delta,\alpha)_n$ set $E^*$ (with the implicit constants in the definition of a 
$(\delta,\alpha)_n$ set depending on $K$) such that
$$ E \subseteq E^* \cup \bigcup_{\delta < \delta' \leq 1} E_{\delta'}.$$
\end{Refinement}

\begin{proof}
Define the sets $E_{\delta'}$ by
$$ E_{\delta'} := \{ x \in \R^n: |E \cap \B(x,\delta')| \geq \delta^{-K\eps} \delta^n 
(\delta'/\delta)^\alpha \}$$
and $E^*$ by
$$ E^* := (E \backslash \bigcup_{\delta < \delta' \leq 1} E_{\delta'}) + \B(0,\delta).$$
The required properties on $E_{\delta'}$ and $E^*$ are then easily verified.
\end{proof}

\begin{Separation}\label{sep-lemma}  Let $X$ be a $(\delta,\alpha)_n$ set in $\R^n$ for some 
$0 < \alpha < n$ such that $|X| \approx \delta^{n-\alpha}$.  Then there exist refinements 
$X_1$, $X_2$ of $X$ which respectively live in cubes $Q_1$, $Q_2$ of size and separation 
$\approx 1$ with $|Q_1| = 
|Q_2|$, and $|X_1|, |X_2| \approx \delta^{n-\alpha}$.
\end{Separation}

\begin{proof} By \eqref{dim} we see that
$$ |X \cap Q| \leq 10^{-n} |X|$$
for all cubes $Q$ of side-length $\delta^{\Cc \eps}$, if $\Cc$ is a sufficiently large 
constant.  The claim then follows by covering $B(0,C)$ with such cubes, extracting the top 
$5^n$ cubes in that collection which maximize $|X \cap Q|$, picking two of those cubes 
$Q_1$, $Q_2$ which are not adjacent, and setting $X_i := X \cap Q_i$ for $i=1,2$.  We leave
 the verification of the desired properties to the reader.
\end{proof}

For any function $f$ in $\R^2$, define the Kakeya maximal function $f^*_\delta(\omega)$ for 
$\omega \in S^1$ by
$$ f^*_\delta(\omega) := \sup_R \frac{1}{|R|} \int_R |f|,$$
where $R$ ranges over all $1 \times \delta$ rectangles oriented in the direction $\omega$.  

The following estimate can be found in 
\cite{cordoba} (see also Lemma \ref{pseudo-cordoba}):

\begin{Kakeya}[C\'ordoba's estimate]\label{cordoba}
We have
$$ \| f^*_\delta \|_2 \lessapprox \|f\|_2.$$
Dually, if we set $R_\omega$ be a collection of $\delta \times 1$ rectangles oriented in a 
$\delta$-separated set of directions, then 
$$ \| \sum_{\omega} \chi_{R_\omega} \|_2 \lessapprox 1.$$
\end{Kakeya}

\section{Arithmetic combinatorics}\label{arithmetic}

We shall prove Theorem \ref{grand-equiv} by showing that
\begin{align*}
\hbox{Bilinear Distance} &\implies \hbox{Discretized Ring } 
\implies \\
\hbox{ Discretized Furstenburg } &\implies \hbox{Bilinear Distance }.
\end{align*}

We shall need a number of standard results concerning the cardinality of sum-sets $A+B$ and 
difference sets $A-B$, and partial sum-sets $\{ a+b: (a,b) \in G \}$, where $G$ is a large 
subset of $A \times B$.

We first give the results in a discrete setting.

\begin{lemma}\cite{ruzsa} Suppose $A_1, A_2$ are finite subsets of $\R$ such that
$$ \#(A_1 + A_2) \approx \# A_1 \approx \# A_2.$$
Then we have
$$ \#(A_{i_1} \pm \ldots \pm A_{i_N}) \approx \# A_1$$
for all choices of signs $\pm$ and $i_1, \ldots, i_N \in \{1,2\}$, where the implicit 
constants depend on $N$.  Also, we can find a refinement $A'_1$ of $A_1$ and a real number 
$x$ such that $x+A'_1$ is a refinement of $A_2$.
\end{lemma}

\begin{proof} Most of these results are in \cite{ruzsa}.  For the last result, observe that 
the discrete function $\chi_{-A_1} * \chi_{A_2}$ has an $l^1$ norm $\approx (\# A_1)^2$ and 
is supported in a set of cardinality $\approx \# A_1$ by the results in \cite{ruzsa}.  Thus 
one can find an $x$ such that $\chi_{-A_1} * \chi_{A_2}(x) \gtrapprox \# A_1$, and the claim 
follows by setting $A'_1 = A_1 \cap (A_2-x)$.
\end{proof}

We also need Bourgain's variant of the Balog-Szemer\'edi theorem \cite{borg:high-dim} 
(as used in Gowers \cite{gowers:isgod} ), namely

\begin{lemma}\cite{borg:high-dim} Let $N \gg 1$ be an integer, and let $A$, 
$B$ be finite subsets of $\R$ such that
$$ \# A, \# B \approx N.$$
Suppose there exists a refinement $G$ of $A \times B$ such that
$$ \# \{ a+b: (a,b) \in G \} \lessapprox N.$$
Then we can find refinements $A'$, $B'$ of $A$ and $B$ respectively
such that $G \cap (A' \times B')$ is a refinement of $A' \times B'$, and for all 
$(a',b') \in A' \times B'$ we have
$$ \# \{ (a_1,a_2,a_3,b_1,b_2,b_3) \in A \times A \times A \times B \times B \times B: 
a'-b' = (a_1-b_1) - (a_2 - b_2) + (a_3 - b_3) \} \approx N^5.$$
In particular, we have
$$ \# (A' - B') \approx N.$$
\end{lemma}

We can easily replace these discrete lemmata with $\delta$-discretized variants as follows.

\begin{corollary}\label{ruzsa-lemma} Suppose $A, B$ are finite unions of intervals of 
length
 $\approx \delta$ such that
$$ |A+B| \approx |A| \approx |B|.$$
Then we have
$$ |A \pm \ldots \pm A| \approx |A|$$
for all choices of signs $\pm$, with the implicit constants depending on the number 
of signs.  Also, we can find a refinement $A'$ of $A$ and a real number $x$ such that
 $x+A'$ is a refinement of $B$.
\end{corollary}

\begin{Perfection}\label{borg-lemma} Let $r \gg \delta$, and let $A$, $B$ be finite unions 
of intervals of length $\approx \delta$ such that
$$ |A|, |B| \approx r.$$
Suppose there exists a refinement $G$ of $A \times B$ such that
$$ |\{ a+b: (a,b) \in G \}| \lessapprox N.$$
Then we can find $\delta$-discretized refinements $A'$, $B'$ of $A$ and $B$ respectively
such that $G \cap (A' \times B')$ is a refinement of $A' \times B'$, and for all 
$(a',b') \in A' \times B'$ we have
$$ |\{ (a_1,a_2,a_3,b_1,b_2,b_3) \in A \times A \times A \times B \times B 
\times B: a'-b' = (a_1-b_1) - (a_2 - b_2) + (a_3 - b_3) \}| \approx r^5.$$
In particular, we have
$$ |A' - B'| \approx r.$$
\end{Perfection}

To obtain these corollaries, we first observe that any $\delta$-discretized set $A$ contains 
the $\approx \delta$-neighbourhood of a discrete set $A^*$ of cardinality 
$\# A^* \approx |A|/\delta$ which is contained in an arithmetic progression of spacing 
$\approx \delta$.  The claims then follow by applying the previous lemmata to 
$A^*$, $B^*$. (See also the proof of \cite{borg:high-dim}, Lemma 2.83).

We also observe the trivial estimate
\be{triv-sum}
|A + B| \gtrapprox |A|, |B|
\end{equation}
for all sets $A$, $B$.

If we also assume that the sets $A$, $B$ are contained in the annulus $\A$ then one can 
also obtain analogues of \eqref{triv-sum} and the above two Corollaries in which addition 
and subtraction are replaced by multiplication and division respectively.  This simply 
follows by applying a logarithmic change of variables.  In the next section we shall use the 
fact that multiplication distributes over addition, to obtain hybrid versions of the above 
results. 

\section{Bilinear Distance Conjecture \ref{less-naive} implies Ring Conjecture \ref{erdos}}

Assume that the Bilinear Distance Conjecture \ref{less-naive} is true for some 
absolute constant $c_1 > 0$.  In 
this section we show how the Ring Conjecture \ref{erdos} follows.

\begin{figure}[htbp]\label{fig3} \centering
\ \psfig{figure=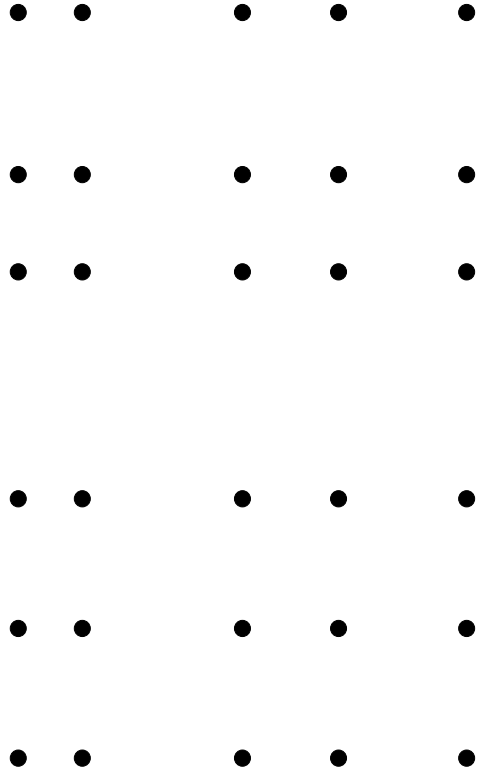,height=4in}
\caption{A set which contradicts the distance conjecture if a half-dimensional
ring exists and constitutes its vertical and horizontal sets of projections.}
\end{figure}

Let $0 <  \eps \ll 1$ be fixed.  We may assume that $\delta$ is sufficiently small depending 
on $\eps$, since the Ring Conjecture is trivial otherwise.   We may also assume that 
$\delta$ 
is dyadic.  Assume for contradiction that one can find a $(\delta,\frac{1}{2})_1$-set $A 
\subset \A$ of measure $|A| \approx \delta^{1/2}$ such that
\be{amult}
|A + A|, |A\cdot A| \lessapprox \delta^{1/2}
\end{equation}
We will obtain a contradiction from this, and it will be clear from the nature of the 
argument that one can in fact show that at least one of $A+A$, $A \cdot A$ has measure 
$\gtrapprox \delta^{\frac{1}{2}-c_4}$ for some absolute constant $c_4 > 0$ depending on 
$c_1$.

From Separation \ref{sep-lemma} one can find refinements $A_1$, $A_2$ of $A$ which are 
contained 
in intervals of size and separation $\approx 1$ and have measure $|A_1|, |A_2| \approx 
\delta^{1/2}$. From the additive and multiplicative versions of \eqref{triv-sum} we thus 
have
\be{init}
 |A_1|, |A_2|, |A_1 + A_2|, |A_1 A_2| \approx \delta^{1/2}.
\end{equation}

Heuristically, the idea is to apply the Bilinear Distance Conjecture \ref{less-naive} with 
$E_0$, $E_1$, $E_2$ 
equal to $A_1 \times A_1$, $A_1 \times A_2$, $A_2 \times A_1$ respectively.  The difficulty 
with this is that we cannot quite control the distance set $\sqrt{(A_1-A_1)^2 + 
(A_1-A_2)^2}$ accurately using \eqref{init}, however this difficulty can be avoided if we 
pass to various refinements of $A$.

We turn to the details.  
From \eqref{init} and Perfection \ref{borg-lemma} with $A,B,G$ set to $A_1$, $A_2$, $A_1 
\times A_2$ respectively, and some re-labeling, we can find $\delta$-discretized 
refinements $C$, $D$ of $A_1$, $A_2$ respectively such that
\be{sum}
|\{ (a_1,a_2,a_3,a_4,a_5,a_6) \in A^{\oplus 6}: d-c = (a_1 - a_4) - (a_2 - a_5) +
 (a_3 - a_6) \}| \approx \delta^{5/2}
\end{equation}
for all $(c,d) \in C \times D$.  From construction we have
\be{sep}
|c-d| \approx 1 \hbox{ for all } c \in C, d \in D.
\end{equation}

\begin{lemma}\label{acd}  We have
$$ |A \cdot A \cdot A \cdot (C-D) / (A \cdot A)| = | \{ 
\frac{a_1 a_2 a_3}{a_4 a_5}(c-d): a_1,a_2,a_3,a_4,a_5 \in A, c \in C, d \in D \}| 
\approx \delta^{1/2}.$$
\end{lemma}

\begin{proof}
The lower bound is clear from \eqref{sep} and the multiplicative version of 
\eqref{triv-sum}, so it suffices to show the upper bound.

Fix $a_1,a_2,a_3,a_4,a_5$, $c$, $d$.  By multiplying \eqref{sum} by 
$a_1 a_2 a_3/a_4 a_5$, which is $\approx 1$, we see that
$$
|\{ (e_1,e_2,e_3,e_4,e_5,e_6) \in (A \cdot A \cdot A \cdot A / (A \cdot A))^{\oplus 6}:
 \frac{a_1 a_2 a_3}{a_4 a_5}(d-c) = (e_1 - e_4) - (e_2 - e_5) + (e_3 - e_6) \}| \gtrapprox 
\delta^{5/2}.
$$
Integrating this over all possible values of $\frac{a_1 a_2 a_3}{a_4 a_5} (d-c)$ and using 
Fubini's theorem we obtain
$$ |A \cdot A \cdot A \cdot A / (A \cdot A)|^6 \gtrapprox \delta^{5/2} |A \cdot (C-D)|.$$
On the other hand, from \eqref{amult} and the multiplicative form of Corollary 
\ref{ruzsa-lemma} we have
$$ | A \cdot A \cdot A \cdot A / (A \cdot A) | \approx \delta^{1/2}.$$
The claim follows by combining the above two estimates.
\end{proof}

From \eqref{init} and the multiplicative version of \eqref{triv-sum} we have
$$ 
|C|, |D|, |CD| \approx \delta^{1/2}.$$
From the multiplicative form of Perfection \ref{borg-lemma} with $A := C$ and $B:=1/D$, we 
may thus find refinements $C'$, $D'$ of $C$, $D$ respectively such that
\be{mult}
|\{ (c_1,c_2,c_3,d_1,d_2,d_3) \in C \times C \times C \times D \times D \times D: 
cd = (c_1d_1) (c_2 d_2)^{-1} (c_3 d_3) \}| \approx \delta^{5/2}
\end{equation}
for all $c \in C'$, $d \in D'$.  

\begin{lemma}\label{cdcd} We have
$$
|C'D'-C'D'| = | \{ cd - c'd' : c,c' \in C', d,d' \in D' \} | \approx \delta^{1/2}.
$$
\end{lemma}

\begin{proof}
As before, the lower bound is immediate from the additive and multiplicative versions of 
\eqref{triv-sum}, so it suffices to show the upper bound.

Fix $c,c',d,d'$, and let $X$ denote the set in \eqref{mult}.
Observe that for all $(c_1,c_2,c_3,d_1,d_2,d_3) \in X$ we have the telescoping identity
$$ cd - c'd' = x_1 - x_2 + x_3 - x_4$$
where
\begin{align*}
x_1 &:= \frac{ (c_1 - d') d_1 c_3 d_3 }{c_2 d_2}\\ 
x_2 &:= \frac{ d' (c' - d_1) c_3 d_3 }{c_2 d_2} \\
x_3 &:= \frac{ d' c' (c_3 - d_2) d_3 }{c_2 d_2} \\
x_4 &:= \frac{ d' c' d_2 (c_2 - d_3) }{c_2 d_2}.
\end{align*}
Indeed, we have the identities
\begin{align*}
\frac{c_1 d_1 c_3 d_3}{c_2 d_2} &= cd \\
\frac{d' d_1 c_3 d_3}{c_2 d_2} &= cd - x_1\\
\frac{d' c' c_3 d_3}{c_2 d_2} &= cd - x_1 + x_2\\
\frac{d' c' d_2 d_3}{c_2 d_2} &= cd - x_1 + x_2 - x_3\\
c'd' = \frac{d' c' d_2 c_2}{c_2 d_2} &= cd - x_1 + x_2 - x_3 + x_4.
\end{align*}
As a consequence of these identities, \eqref{sep} and some algebra we see the
$$ (c_1, c_2, c_3, d_1, d_2, d_3) \mapsto (x_1, x_2, x_3, x_4, c_2, d_2)$$
is a diffeomorphism on $X$ (recall that $c$, $d$, $c'$, $d'$ are fixed).  From \eqref{mult} 
we thus have
$$ | \{ (x_1,x_2,x_3,x_4,c_2,d_2) \in (A \cdot A \cdot A \cdot (C-D)/(A \cdot A))^{\oplus 4} 
\times C \times D: cd - c'd' = x_1 - x_2 + x_3 - x_4 \} | \gtrapprox \delta^{5/2}.$$
Integrating this over all values of $cd - c'd'$ and using Fubini's theorem we obtain 
$$ |C'D'-C'D'| \gtrapprox \delta^{5/2} | A \cdot A \cdot A \cdot (C-D)/(A \cdot A) |^4 |C| 
|D|.$$
The claim then follows from Lemma \ref{acd}.
\end{proof}

From the above lemma and the multiplicative form of \eqref{triv-sum} we have
$$ |C'|, |D'|, |C'D'| \approx \delta^{1/2}.$$
From the multiplicative version of Corollary \ref{ruzsa-lemma} we can therefore find a 
refinement $F$ of $C'$ and a real number $x \approx 1$
such that $xF$ is a refinement of $D'$.  In particular, since $FF - FF$ is a subset of 
$x^{-1}(C'D'-C'D')$, we thus see that
$$ |FF - FF| \approx |FF| \approx |F| \approx \delta^{1/2}.$$
From Corollary \ref{ruzsa-lemma} we thus have
$$ |FF - FF - FF + FF + FF - FF - FF + FF| \approx \delta^{1/2}.$$
Since $(F-F)^2 \subset FF - FF - FF + FF$, we thus have
$$ |(F-F)^2 + (F-F)^2| \lessapprox \delta^{1/2}.$$
The set $F$ is a $(\delta,\frac{1}{2})_1$ set with measure $\approx \delta^{1/2}$.  From 
Separation \ref{sep-lemma} we may find refinements $F_1, F_2$ of $F$ which are contained in 
intervals $I_1$, $I_2$ of size and separation $\approx 1$ such that $|I_1| = |I_2|$ and 
$|F_1|, |F_2| \approx \delta^{1/2}$.  

Define
$$E_0 := F_1 \times F_1, E_1 := F_1 \times F_2, E_2 := F_2 \times F_1, \quad
Q_0 := I_1 \times I_1, Q_1 := I_1 \times I_2, Q_2 := I_2 \times I_1.$$
It is clear that $Q_0, Q_1, Q_2$ obey \eqref{gen-sep} and that $E_0, E_1, E_2$ are 
$(\delta,1)_2$ sets of measure $\approx \delta$ contained in $Q_0$, $Q_1$, $Q_2$ 
respectively.

Let $D$ denote the set
$$ D = \sqrt{(F_2 - F_1)^2 + (F_1 - F_1)^2}.$$
Clearly $D$ is a $\delta$-discretized set of measure $|D| \lessapprox \delta^{1/2}$ which 
lives in $\A$.  In fact, from the size and separation of $F_1$ and $F_2$ we have 
\be{dsize} |D| \approx \delta^{1/2}.\end{equation}
Also, we have
$$ |x_1 - x_0|, |x_2 - x_0| \in D$$
for all $x_0 \in E_0$, $x_1 \in E_1$, $x_2 \in E_2$.  In particular, we have
\be{skim}
|\{ (x_0,x_1,x_2) \in E_0 \times E_1 \times E_2: |x_0-x_1|, |x_0-x_2| \in D \}| = |E_0| 
|E_1| |E_2| \approx \delta^3.
\end{equation}
We are almost ready to apply the hypothesis \eqref{eps1}, however the one thing which is 
missing is that $D$ need not satisfy \eqref{dim}.  To rectify this we shall remove some 
portions from $D$.

Apply Refinement \ref{refine} to obtain a covering
$$ D \subset D^* \cup \bigcup_{\delta < \delta' \ll 1} D_{\delta'}$$
with the properties asserted in Refinement \ref{refine}  , 
and $\Ca$ equal to a large constant to be chosen 
shortly.

\begin{proposition}  For all $\delta' > \delta$, we have
$$ |\{ (x_0,x_1) \in E_0 \times E_1: |x_0 - x_1| \in D_{\delta'} \}| \lessapprox \delta^2 
\delta^{\Ca\eps/100}$$
and
$$ |\{ (x_0,x_2) \in E_0 \times E_2: |x_0 - x_2| \in D_{\delta'} \}| \lessapprox \delta^2 
\delta^{\Ca\eps/100}.$$
\end{proposition}

\begin{proof} Fix $\delta'$.  We may assume that $\eps$ is sufficiently small depending on 
$\Ca$, and $\delta$ is sufficiently small depending on $\Ca$ and $\eps$, since the claim is 
trivial otherwise.

By reflection symmetry it suffices to prove the first estimate.  Suppose for contradiction 
that
$$
|\{ (x_0,x_1) \in E_0 \times E_1: |x_0 - x_1| \in D_{\delta'} \}| \gtrapprox \delta^2 
\delta^{\Ca\eps/100}.$$
From Cauchy-Scwartz \ref{cauchy} we thus have
$$
|\{ (x_0,x_1,x'_1) \in E_0 \times E_1 \times E_1: |x_0 - x_1| \in D_{\delta'}, 
|x_0 - x'_1| \in D_{\delta'} \}| \gtrapprox \delta^3 \delta^{\Ca\eps/50}.$$
Write $x_1 = (\x_1,\y_1)$, $x'_1 = (\x'_1,\y'_1)$.  Observe that
$$
|\{ (x_0,x_1,x'_1) \in E_0 \times E_1 \times E_1: |\x_1 - \x'_1| \lessapprox 
\delta^{\Ca \eps/10} \}| \lessapprox \delta^3 \delta^{\Ca \eps/20}.$$
This is because for fixed $\x_1$, $\x'_1$ can only range in a set of measure 
$\lessapprox \delta^{1/2} \delta^{\Ca \eps/20}$ thanks to \eqref{dim} and the fact that 
$F_1$ is a $(\delta, \frac{1}{2})_1$ set.  Subtracting the two inequalities we obtain 
(if $\delta$ is sufficiently small)
$$
|\{ (x_0,x_1,x'_1) \in E_0 \times E_1 \times E_1: |x_0 - x_1| \in D_{\delta'}, 
|x_0 - x'_1| \in D_{\delta'}, |\x_1 - \x'_1| \gtrapprox \delta^{\Ca \eps/10} \}| 
\gtrapprox \delta^3 \delta^{\Ca\eps/50}.$$
Since $|E_1| \approx \delta$, we may thus find $x_1, x'_1 \in E_1$ such that \be{x-sep}
|\x_1 - \x'_1| \gtrapprox \delta^{\Ca \eps/10}
\end{equation}
and
\be{x0}
|\{ x_0 \in E_0: |x_0 - x_1| \in D_k, |x_0 - x'_1| \in D_{\delta'} \}| \gtrapprox \delta 
\delta^{\Ca\eps/50}.
\end{equation}
From Refinement \ref{refine} $D_{\delta'}$ can be covered by $\lessapprox \delta^{\Ca\eps} 
{\delta'}^{-1/2} $ intervals in $\A$ of length $\lessapprox \delta'$.  From this fact, 
\eqref{x-sep}, and the geometry of annuli which intersect non-tangentially,
we see that the set in \eqref{x0} can be covered by
$\lessapprox {\delta'}^{-1} \delta^{2\Ca\eps}$ balls of radius $\lessapprox \delta^{-\Ca 
\eps/5} \delta' $.  Since $E_0$ is a $(\delta,1)_2$ set, we see from \eqref{dim} that
$$ \hbox{LHS of \eqref{x0}} \lesssim {\delta'}^{-1} \delta^{2\Ca\eps}
\delta' \delta^{-\Ca \eps/5}.$$
But this contradicts \eqref{x0} if $\delta$ is sufficiently small.  This concludes the 
proof of the proposition.
\end{proof}

From \eqref{skim} and the above proposition we see that (if $\Ca$ is a large enough absolute 
constant, and $\delta$ is sufficiently small depending on $\eps$, $\Ca$)
\be{d-b}
|\{ (x_0,x_1,x_2) \in E_0 \times E_1 \times E_2: |x_0-x_1|, |x_0-x_2| \in D^* \}|  
\gtrapprox \delta^3.
\end{equation}
From \eqref{dsize} we have $|D^*| \lessapprox \delta^{1/2}$.  From elementary geometry and a 
change of variables we have
$$ |\{ x_0 \in E_0: |x_0-x_1|, |x_0-x_2| \in D^* \}| \lessapprox |D^*|^2$$
for all $x_1 \in E_1$, $x_2 \in E_2$.  Integrating this over $x_1$ and $x_2$ and comparing 
with the previous we thus see that $|D^*| \approx \delta^{1/2}$. But then \eqref{d-b} 
contradicts \eqref{eps1} (with $D$ replaced by $D^*$), if $\eps$ is sufficiently small 
depending on $c_1$ and $\delta$ sufficiently small depending on $\eps$.  The full claim of 
the proposition follows by a modification of this argument, providing that $c_4$ is 
sufficiently small depending on $c_1$.

\section{Ring Conjecture \ref{erdos} implies Discretized Furstenburg 
Conjecture \ref{furst-disc}}\label{erdos-furst}

Assume that the Ring Conjecture \ref{erdos} is true for some absolute constant $c_4 > 0$.  In this 
section we show how the Discretized Furstenburg Conjecture \ref{furst-disc} follows.

The main idea is that $R$ is a half-dimensional ring then $R \times R$ contains
a one dimensional set of lines each of which contain half dimensional sets.
That many of these lines are parallel seems hardly consequential and we
will deal with it by an appropriately chosen projective transformation.

Let $0 <  \eps \ll 1$ be fixed.  We may assume that $\delta$ is sufficiently small depending 
on $\eps$, since \eqref{furst-est} is trivial otherwise, and may assume $\delta$ is dyadic 
as before.   Let $E$, $\Omega$, $R_\omega$ be as in the Discretized Furstenburg
Conjecture \ref{furst-disc}.  Assume for 
contradiction that 
\be{hyp} | \{ (x_0,x_1) \in E \times E: x_1, x_0 \in R_{\omega} \hbox{ for some } \omega 
\in \Omega \} |
\gtrapprox \delta^2
\end{equation}
We will obtain a contradiction from this, and it will be clear from the nature of the 
argument that \eqref{furst-est} in fact holds for some absolute constant $c_3 > 0$ depending 
on $c_4$.

It will be convenient to define the non-transitive relation $\sim$ by defining $x \sim y$ 
if and only if $x, y \in R_\omega$ for some $\omega \in \Omega$.  We also write $x_1, 
\ldots, x_n \sim y_1, \ldots y_m$ if $x_i \sim y_j$ for all $1 \leq i \leq n$ and all 
$1 \leq j \leq m$.

From \eqref{hyp} we then have
\be{hyp2} | \{ (x_0, x_1) \in E \times E: x_0 \sim x_1 \} | \gtrapprox \delta^2.
\end{equation}
 
Roughly speaking, the idea will be to find $x_1, x'_1 \in E$ and a refinement $E''$ of $E$ 
such that $x_0 \sim x_1, x_0 \sim x'_1$ for all $x_0 \in E''$, and such that there are many 
relations between pairs of points in $E''$.  Then after a projective transformation sending 
$x_1$, $x'_1$ to the cardinal points at infinity we can transform $E''$ to a Cartesian 
product of two $(\delta,\frac{1}{2})_1$ sets of measure $\approx \delta^{1/2}$, at which 
point the ring structure of these sets can be easily extracted.

We turn to the details.  From \eqref{hyp2} and the fact that $|E| \approx \delta$, we see 
that
\be{hyp3} | \{ (x_0, x_1) \in E' \times E: x_0 \sim x_1 \} | \gtrapprox \delta^2
\end{equation}
where
$$ E' = \{ x_0: | \{ x_1 \in E: x_0 \sim x_1 \}| \approx \delta \}$$
provided the constants are chosen appropriately.

Let $\Cb$ be a large constant to be chosen later, and let $E_1$ be the set 
$$ E_1 = \{ x_1 \in E: \sum_{\omega \in \Omega} \chi_{R_\omega}(x_1) \leq \delta^{-\Cb \eps} 
\delta^{-1/2} \}.$$

From Kakeya \ref{cordoba} and Chebyshev we have
$$ |E \backslash E_1| \lesssim \delta^{2\Cb \eps} \delta$$
and thus
$$ | \{ (x_0, x_1) \in E' \times (E \backslash E_1): x_0 \sim x_1 \} | \lessapprox 
\delta^{2\Cb \eps} \delta^2.
$$

If we then choose $\Cb$ is large enough, and $\delta$ is small enough depending on $\Cb$ 
and $\eps$, we thus see from \eqref{hyp3} that
\be{hyp4} | \{ (x_0, x_1) \in E' \times E_1: x_0 \sim x_1 \} | \gtrapprox \delta^2.
\end{equation}
In particular, we have $|E_1| \approx \delta$ as before.  Henceforth $\Cb$ is fixed so that 
\eqref{hyp4} applies.

From \eqref{hyp4} and Cauchy-Schwarz \ref{cauchy} we have
\be{hyp5} | \{ (x_0, x_1,x'_1) \in E' \times E_1 \times E_1: x_0 \sim x_1,x'_1 \} | 
\gtrapprox \delta^3.
\end{equation}

Let $\Cd$ be a large constant to be chosen later.

\begin{figure}[htbp] \centering
\ \psfig{figure=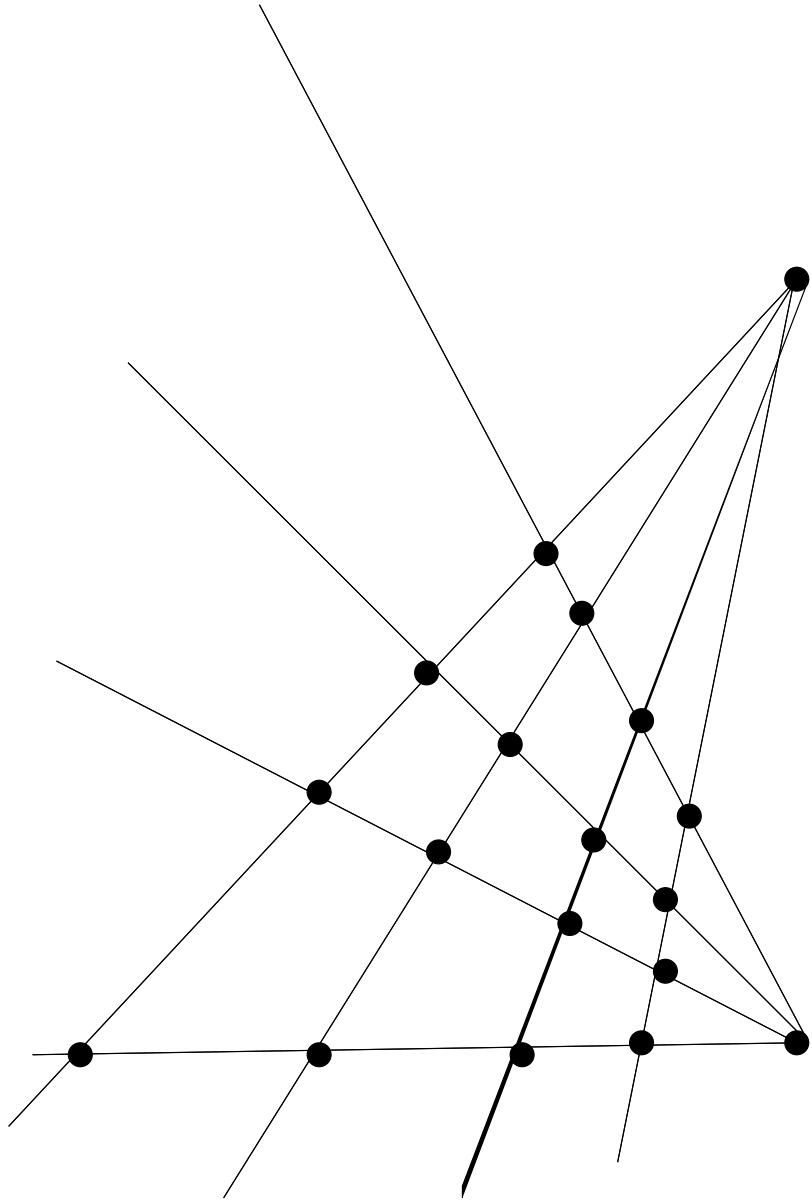}
\caption{A Furstenburg set when viewed from $x_1$ and $x_1^{\prime}$.  Note how
this resembles a projective transformation of Figure \ref{fig3}. }
\end{figure}

\begin{lemma} If $\Cd$ is large enough, and $\delta$ is small enough depending on $\Cd$ and 
$\eps$, we have 
$$ 
|\{ (x_0, x_1, x'_1) \in E' \times E_1 \times E_1: x_0 \sim x_1,x'_1; |(x_1 - x_0) \wedge 
(x'_1 - x_0)| \geq \delta^{\Cd \eps} \}| \approx \delta^3.$$
\end{lemma}

\begin{proof}
From \eqref{hyp4} it suffices to show that
\be{hyp4-neg}
|\{ (x_0, x_1, x'_1) \in E' \times E_1 \times E_1: x_0 \sim x_1,x'_1; |(x_1 - x_0) \wedge 
(x'_1 - x_0)| \leq \delta^{\Cd \eps} \}| \lessapprox \delta^{\Cd \eps/8} \delta^3.
\end{equation}
(The constant $8$ is non-optimal, but this is irrelevant for our purposes).  
In order to have
$$ |(x_1 - x_0) \wedge (x'_1 - x_0)| \leq \delta^{\Cd \eps}$$
one must either have $|x_1 - x'_1| \lessapprox \delta^{\Cd \eps/2}$, or
that $|x_1 - x'_1| \gtrapprox \delta^{\Cd \eps/2}$ and $x_0$ is within $\lessapprox 
\delta^{\Cd \eps/2}$ of the line joining $x_1 - x'_1$.

Let us consider the contribution of the former case.  Since $E_1$ is a $(\delta,1)_2$ set, 
we see that each pair $(x_0,x_1)$ contributes a set of measure $\lessapprox 
\delta^{\Cd \eps/2} \delta$ to \eqref{hyp4-neg}.  From Fubini's theorem we thus see that the
 contribution of this case to \eqref{hyp4-neg} is acceptable.

Now let us consider the contribution of the latter case.  By Fubini's theorem again it 
suffices to show that
$$
|\{ x_0 \in E': x_0 \in S \}| \lessapprox \delta^{\Cd \eps/8} \delta$$
for any strip $S$ of width $\lessapprox \delta^{\Cd \eps/2}$.

Fix $S$.  From the definition of $E'$ and Fubini's theorem it suffices to show that
$$
|\{ (x_0, x_2) \in E' \times E: x_0 \in S, x_2 \sim x_0 \}| \lessapprox \delta^{\Cd \eps/8} 
\delta^2.
$$
From the definition of $\sim$, we can estimate the left-hand side by
$$ \sum_{\omega \in \Omega} |S \cap R_\omega| |R_\omega|.$$
Since $R_\omega$ is a $(\delta,\frac{1}{2})_1$ set, we have $|R_\omega| \lessapprox 
\delta^{3/2}$.  Also, if $\omega$ makes an angle of $\gtrapprox \delta^{\Cd \eps/4}$ with 
$S$ we have $|S \cap R_\omega| \lessapprox \delta^{\Cd \eps/8} \delta^{3/2}$ by \eqref{dim}
 and elementary geometry, otherwise we may estimate $|S \cap R_\omega| \leq |R_\omega| 
\lessapprox \delta$.  Inserting these estimates into the previous and using the 
$\delta$-separated nature of the $\omega$, we see that
$$ \sum_{\omega \in \Omega} |S \cap R_\omega| |R_\omega| \lessapprox \delta^{\Cd \eps/8} 
\delta^2$$
as desired.
\end{proof}

Henceforth $\Cd$ is fixed so that the above lemma applies.  We now suppress all explicit 
mention of $\Cb$, $\Cd$ and absorb these factors into the $\lessapprox$ notation.

From the lemma and the fact that $|E_1| \approx \delta$, we can thus find $x_1, x'_1$ in 
$E_1$ such that
$$ |\{ x_0 \in E': x_0 \sim x_1,x'_1; |(x_1 - x_0) \wedge (x'_1 - x_0)| \approx 1 \}| 
\gtrapprox \delta.$$

Fix $x_1$, $x'_1$.  Clearly one must have $|x_1 - x'_1| \approx 1$, else the left-hand side 
is necessarily zero.  If we define $Q$ by
$$ Q := \{ x_0 \in \R^2: |x_0| \lessapprox 1; |(x_1 - x_0) \wedge (x'_1 - x_0)| \approx 1 \}$$
and $E''$ by $E'' := E' \cap Q$,
then clearly $|E''| \approx \delta$ if
we have chosen the origin appropriately.  Also, if we define $\Omega_1$, $\Omega'_1$ by
$$ \Omega_1 := \{ \omega \in \Omega: x_1 \in R_\omega \}, \quad \Omega'_1 := 
\{ \omega \in \Omega: x'_1 \in R_\omega \},$$
then we have 
$$ \sum_{\omega \in \Omega_1} \sum_{\omega' \in \Omega'_1}
|E'' \cap R_\omega \cap R_{\omega'}| \gtrapprox \delta.$$
From the definition of $E_1$ we note that
$$\# \Omega_1, \# \Omega'_1 \lessapprox \delta^{-1/2}.$$

From the definition of $E'$ we see that
$$\# \{ \omega \in \Omega: x_0 \in R_\omega \cap Q \} \gtrapprox \delta^{-1/2}$$
for all $x_0 \in E''$.  Integrating this on $E''$, which has measure $\approx \delta$, 
we obtain
\be{om}
\sum_{\omega \in \Omega} |R_\omega \cap Q \cap E''| \gtrapprox \delta^{1/2}.
\end{equation}

Let $\Omega_2$ denote those $\omega \in \Omega$ for which the direction of the bounding 
rectangle for $R_\omega$ stays at a distance $\approx 1$ from $x_1$ and $x'_1$.

\begin{lemma}  If the constants in the definition of $\Omega_2$ are chosen appropriately, 
we have
\be{om2}
\sum_{\omega \in \Omega_2} |R_\omega \cap Q \cap E''| \gtrapprox \delta^{1/2}.
\end{equation}
\end{lemma}

\begin{proof}
From elementary geometry we see that
$$ \chi_{E''}^*(\omega') \gtrapprox \delta^{-1} |R_\omega \cap Q \cap E''|$$
whenever $|\omega' - \omega| \ll \delta$, where $\chi_{E''}^*$ is the Kakeya maximal 
function of $\chi_{E''}$.  From Kakeya \ref{cordoba} we thus have
$$ \delta \sum_{\omega \in \Omega} \delta^{-2} |R_\omega \cap Q \cap E''|^2
\lessapprox |E''| \approx \delta.$$
From Cauchy-Schwarz we thus have
$$ \sum_{\omega \in \Omega \backslash \Omega_2} |R_\omega \cap Q \cap E''|
\lessapprox \# (\Omega \backslash \Omega_2)^{1/2} \delta.$$
If one defines the constants in $\Omega_2$ appropriately, the claim then follows from 
\eqref{om}.
\end{proof}

Let $\R\P^2$ denote the projective plane, i.e. the points in $\R^3 \backslash \{0\}$ with 
$x$ identified with $tx$ for all $t \in \R \backslash \{0\}$.  We embed $\R^2$ into 
$\R\P^2$ in usual manner, identifying $(x,y)$ with $(x,y,1)$.

Let $L: \R\P^2 \to \R\P^2$ be a projective linear transformation which sends $x_1$ to 
$(1,0,0)$ and $x'_1$ to $(0,1,0)$, but maps $Q$ to a subset of $\B(0,1)$ with Jacobian 
$\approx 1$ on $Q$.  (This is possible because of the construction of $Q$).  In particular 
we have 
\be{led}
|L(E'')| \approx \delta.
\end{equation}

The set $\bigcup_{\omega \in \Omega_1} R_\omega \cap Q$ stays a distance $\approx 1$ from 
the line joining $x_1$ and $x_2$, and is also contained in the union of $\lessapprox 
\delta^{-1/2}$ rectangles of dimensions about $\delta \times 1$ which pass through $x_1$.  
From this fact and some elementary projective geometry we see that
$$ L(\bigcup_{\omega \in \Omega_1} R_\omega \cap Q) \subset \R \times B$$
for some $\delta$-discretized set $B \subset [-1,1]$ with $|B| \approx \delta^{1/2}$.  
Similarly we have
$$ L(\bigcup_{\omega \in \Omega'_1} R_\omega \cap Q) \subset A \times \R$$
for some $\delta$-discretized set $A \subset [-1,1]$ with $|A| \approx \delta^{1/2}$.  
Combining these two facts with the definition of $E''$ we thus have
\be{lab}
 L(E'') \subset A \times B.
\end{equation}

The sets $A$ and $B$ are already our prototypes for half-dimensional rings.
In what follows we refine their geometric properties and establish their
algebraic ones.

For all $\omega \in \Omega_2$, let $\tilde R_\omega$ denote the set
$$ \tilde R_\omega := L(R_\omega \cap Q).$$
From the hypothesis on $R_\omega$ and some elementary projective geometry we see that 
$\tilde R_\omega$ is a $(\delta,\frac{1}{2})_2$ set which is contained in a rectangle of 
dimensions $\approx 1 \times \delta$, and whose long side is oriented at an angle of 
$\approx 1$ to the cardinal directions $(0,\pm 1)$, $(\pm 1,0)$.  From \eqref{om2}, 
\eqref{lab} we have
\be{rab} \sum_{\omega \in \Omega_2} |\tilde R_\omega \cap (A \times B)| \gtrapprox 
\delta^{1/2}.
\end{equation}

The sets $A$ and $B$ need not be $(\delta,\frac{1}{2})_1$ sets because there is no reason 
why they should satisfy \eqref{dim}.  To rectify this we shall refine $A_0$, $B_0$ 
slightly.

Apply Refinement \ref{refine}, with $\Ce$ a large constant to be chosen shortly, to obtain a 
covering
$$ A \subset A^* \cup \bigcup_{\delta < \delta' \leq 1} A_{\delta'}.$$
From Refinement \ref{refine}, $A_{\delta'}$ can be covered by $\lessapprox \delta^{\Ce \eps} 
{\delta'}^{-1/2}$ intervals $I$ of length $\delta'$.  For each such interval $I$ we have
$$ |\tilde R_\omega \cap (I \times B)| \lessapprox {\delta'}^{1/2} \delta;$$
this follows from the properties of $\tilde R_\omega$, \eqref{dim}, and some elementary 
geometry.  Summing this over $I$ and $\omega$, we obtain
$$ \sum_{\omega \in \Omega_2} |\tilde R_\omega \cap (A_{\delta'} \times B)| \lessapprox 
\delta^{\Ce \eps} \delta^{1/2}.$$
Summing this over all $\delta'$, we obtain (if $\Ce$ is sufficiently large, and $\delta$ 
sufficiently small depending on $\Ce$ and $\eps$)
$$ \sum_{\omega \in \Omega_2} |\tilde R_\omega \cap (A^* \times B)| \gtrapprox 
\delta^{1/2}$$
By breaking $A^*$ up into intervals $I$ of length $\approx \delta$ and arguing as before 
we see that
$$ \sum_{\omega \in \Omega_2} |\tilde R_\omega \cap (A^* \times B)| \lessapprox 
|A^*|;$$
since $|A^*| \leq |A| \lessapprox \delta^{1/2}$, we thus see that $|A^*| \approx 
\delta^{1/2}$.  Also, by Lemma \ref{refine} we see that $A^*$ is a 
$(\delta,\frac{1}{2})_1$ set.

By repeating the above argument in the second co-ordinate, we may also find a 
$(\delta,\frac{1}{2})_1$ set $B^*$ of measure $\approx \delta^{1/2}$ such that
$$ \sum_{\omega \in \Omega_2} |\tilde R_\omega \cap (A^* \times B^*)| \approx 
\delta^{1/2}.$$

Let $\Omega'_2$ consist of those $\omega \in \Omega_2$ such that
\be{omp2-def}
|\tilde R_\omega \cap (A^* \times B^*)| \gtrapprox \delta^{3/2}.
\end{equation}
Since $\# \Omega_2 \lessapprox \delta^{-1}$, we thus see that
$$ \sum_{\omega \in \Omega_2 \backslash \Omega'_2} |\tilde R_\omega 
\cap (A^* \times B^*)| \leq \frac{1}{2} \sum_{\omega \in \Omega_2} 
|\tilde R_\omega \cap (A^* \times B^*)|$$
if the constants are chosen correctly.  We thus have
$$ \sum_{\omega \in \Omega'_2} |\tilde R_\omega \cap (A^* \times B^*)| \approx 
\delta^{1/2}.$$
Since $|A^*| \approx \delta^{1/2}$, we can therefore use the pigeonhole principle find 
an $a \in A^*$ such that
$$ \sum_{\omega \in \Omega'_2} |\tilde R_\omega \cap (\{a\} \times B^*)| \gtrapprox 1.$$

Fix such an $a$, and let $\Omega''_2$ consist of those $\omega \in \Omega'_2$ such that
$$ |\tilde R_\omega \cap (\{a\} \times B^*)| \gtrapprox \delta.$$
By repeating the previous argument, we see that
$$ \sum_{\omega \in \Omega''_2} |\tilde R_\omega \cap (\{a\} \times B^*)| \gtrapprox 1$$
for suitable choices of constants.

Let $\Cf$ be a constant to be chosen later.  Since $A^*$ is a $(\delta, \frac{1}{2})_1$ set, 
the set $A^* \cap \B(a,\delta^{\Cf \eps})$ can be covered by $\lessapprox 
\delta^{\Cf \eps/2} \delta^{-1/2}$ intervals $I$ of length $\delta$.  By repeating the 
argument used to refine $A$ and $B$, we have
$$ \sum_{\omega \in \Omega''_2} |\tilde R_\omega \cap ((A^* \cap \B(a,\delta^{\Cf \eps}) 
\times B^*)| \lessapprox \delta^{\Cf \eps/2} \delta^{1/2}.$$
Thus, if $\Cf$ is large enough and $\delta$ is small enough depending on $\Cf$ and 
$\eps$, then
$$ \sum_{\omega \in \Omega''_2} |\tilde R_\omega \cap ((A^* \backslash 
\B(a,\delta^{\Cf \eps}) \times B^*)| \approx \delta^{1/2}.$$
Fix $\Cf$, so that all implicit constants may depend on $\Cf$.  By the pigeonhole principle 
again, one can thus find an $a' \in A^*$ such that $|a-a'| \approx 1$ and
$$ \sum_{\omega \in \Omega''_2} |\tilde R_\omega \cap (\{a'\} \times B^*)| \gtrapprox 1.$$

Fix $a'$.  Let $\Omega'''_2$ consist of those $\omega \in \Omega''_2$ such that
$$ |\tilde R_\omega \cap (\{a'\} \times B^*)| \gtrapprox \delta.$$
Then we have by the same arguments as before that
$$ \sum_{\omega \in \Omega'''_2} |\tilde R_\omega \cap (\{a'\} \times B^*)| \gtrapprox 1.$$
Since $\tilde R_\omega$ is contained in a rectangle of sides $\approx 1 \times \delta$ and 
making an angle of $\approx 1$ with the vertical, we see that
$$ |\tilde R_\omega \cap (\{a'\} \times B^*)| \lessapprox \delta$$
for all $\omega$.  Since $\# \Omega'''_2 \leq \# \Omega \lessapprox \delta^{-1}$, we thus 
see that
\be{om3card} \# \Omega'''_2 \approx \delta^{-1}.
\end{equation}
Consider the set $X \subset B^* \times B^*$ defined by
$$ X := \{ (b,b'): (a,b), (a',b') \in \tilde R_\omega \hbox{ for some } \omega \in 
\Omega'''_2 \}.$$
Each $\omega$ contributes a set of measure $\gtrapprox \delta^2$ to $X$.  Since the 
$\omega$ are $\delta$-separated and $|a-a'| \approx 1$, we see from elementary geometry 
that any given point in $X$ can arise from at most $\lessapprox 1$ values of $\omega$.  
Combining these two facts with \eqref{om3card} we see that
\be{x-bound}
|X| \gtrapprox \delta.
\end{equation}
In particular, $X$ is a refinement of $B^* \times B^*$.

Let $\Cg$ be a large constant to be chosen later.  We now wish to find many values of 
$(b,b') \in X$ and $a'' \in A^*$ such that
\be{lin}
\frac{a''-a'}{a-a'} b + \frac{a-a''}{a-a'} b' \in \tilde B,
\end{equation}
where
$$ \tilde B := B^* + \B(0,C\delta^{1 + \Cg \eps})$$
is a slight enlargement of $B^*$.

\begin{lemma}  If $\Cg$ is a sufficiently large constant, and $\delta$ is sufficiently 
small depending on $\Cg$ and $\eps$, then
$$ |\{ (b,b',a'') \in X \times A^*: |a''-a|, |a''-a'| > \delta^{\Cg \eps}, 
\hbox{\eqref{lin} holds} \}| \gtrapprox \delta^{3/2}.$$
\end{lemma}

\begin{proof}
Fix $(b,b') \in X$.  From \eqref{x-bound} it suffices to show that
$$ |\{ a'' \in A^*: |a''-a|, |a''-a'| > \delta^{\Cg \eps}, \hbox{\eqref{lin} holds} \}| 
\gtrapprox \delta^{1/2}.$$
From the definition of $X$ and the fact that $\Omega'''_2 \subset \Omega'_2$, we can find 
$\omega \in \Omega'_2$ such that $(a,b), (a',b') \in \tilde R_\omega$ and \eqref{omp2-def}
 holds.  From elementary geometry we see $\tilde R_\omega$ stays within $\lessapprox \delta$ 
of the line
$$ \{ (a', \frac{a''-a'}{a-a'} b + \frac{a-a''}{a-a'} b'): a' \in \R \}.$$
Since $\tilde R_\omega$ is $\delta$-discretized, we thus have (if $\Cg$ is sufficiently 
large)
$$ |\{ a'' \in A^*: \hbox{\eqref{lin} holds} \}| \gtrapprox \delta^{1/2}.$$
The separation conditions $|a'' - a|, |a'' - a'| > \delta^{\Cg \eps}$ are easily imposed 
by \eqref{dim}, since $A^*$ is a $(\delta,\frac{1}{2})_1$ set.
\end{proof}

Fix $\Cg$; all implicit constants may now depend on $\Cg$.  Let $T$ denote the set
$$ T =  \{ \frac{a-a''}{a-a'}: a'' \in A^* \} \cap \{ t \in \R: |t| \approx |1-t| 
\approx 1 \};$$
note that $T$ is a $(\delta,\frac{1}{2})_1$ set.  From the above lemma we have
$$ | \{ (b,b',t) \in B^* \times B^* \times T: (1-t)b + tb' \in \tilde B \}| \gtrapprox 
\delta^{3/2}.
$$
Let $B'$ denote the set of all $b' \in B$ such that
\be{juice} |  \{ (b,t) \in B^* \times T: (1-t)b + tb' \in \tilde B \}| \gtrapprox \delta,
\end{equation}
From the previous estimate and the fact that $|B^*| \approx \delta^{1/2}$, we see that $B'$ 
is a refinement of $B^*$ if the constants are chosen correctly.  $B'$ is not quite 
$\delta$-discretized, but this can be easily remedied by introducing the set
$B'' := B' + \B(0,\delta)$.  $B''$ is now a $(\delta,\frac{1}{2})$ set, and every element 
$b' \in B''$ obeys \eqref{juice} if we enlarge the constant $C$ in the definition of 
$\tilde B$ slightly.
In particular, we have
$$ |  \{ (b,b',t) \in B^* \times B'' \times T: (1-t)b + tb' \in \tilde B \}| \gtrapprox 
\delta^{3/2}.
$$
Since $|T| \approx \delta^{1/2}$, there thus exists a $t_0 \in T$ such that
$$ |  \{ (b,b') \in B^* \times B'': (1-t_0)b + t_0 b' \in \tilde B \}| \gtrapprox 
\delta.
$$
Applying Perfection \ref{borg-lemma} with $A$ and $B$ replaced by $(1-t_0)B$ and $t_0 B''$, which 
have measure $\approx \delta^{1/2}$, we can thus find a $\delta$-discretized refinements 
$(1-t_0)B''''$ and $t_0 B'''$ of $(1-t_0)B$ and $t_0 B''$ respectively such that
$$ |(1-t_0)B'''' - t_0 B'''| \approx \delta^{1/2}.$$
From Lemma \ref{ruzsa-lemma} we thus have
$$ |t_0 B''' + t_0 B'''| \approx \delta^{1/2}$$
so that
\be{surfing}
 |B''' + B'''| \approx \delta^{1/2}.
\end{equation}
Note that $B'''$ is a $\delta$-discretized refinement of $B''$ and is therefore a 
$(\delta, \frac{1}{2})_1$ set with measure $\approx \delta^{1/2}$.  

Integrating \eqref{juice} over all $b' \in B'''$ we have
$$ |  \{ (b,b',t) \in B^* \times B''' \times T: (1-t)b + tb' \in \tilde B \}| \gtrapprox 
\delta^{3/2}.
$$
Since $|B^*| \approx \delta^{1/2}$, there thus exists a $b_0 \in B^*$ such that
$$ |  \{ (b',t) \in B''' \times T: (1-t)b_0 + tb' \in \tilde B \}| \gtrapprox \delta.
$$
Fix $b_0$.  We rewrite the above as
$$ |  \{ (f,t) \in (B'''-b_0) \times T: ft \in \tilde B - b_0 \}| \gtrapprox \delta.
$$
Let $\Ch$ be a constant to be chosen later, and define the set
$$ F = \{ f \in B'''-b_0: |f| > \delta^{\Ch \eps} \} + \B(0,\delta).$$
Since $B'''$ is a $(\delta,\frac{1}{2})_1$ set, we have
$$ |(B'''-b_0) \backslash F| \lesssim \delta^{\Ch \eps/2} \delta^{1/2}.$$
In particular, we have $|F| \approx \delta^{1/2}$ and
\be{ft-est}
 |  \{ (f,t) \in F \times T: ft \in \tilde B - b_0 \}| \gtrapprox \delta
\end{equation}
if $\Ch$ is chosen sufficiently large, and $\delta$ sufficiently small depending on $\Ch$ 
and $\eps$. 

Fix $\Ch$; all constants may now depend on $\Ch$.  From the multiplicative form of 
Perfection \ref{borg-lemma} and \eqref{ft-est} we can thus find a $\delta$-discretized 
refinement $F'$ of $F$ such that
$$ |F' F'| \approx \delta^{1/2}.$$
From the previous we also have $|F'+F'| \lessapprox \delta^{1/2}$.  Since $F'$ is a 
$(\delta,\frac{1}{2})_1$ set of measure $\approx \delta^{1/2}$ contained in some annulus 
$\A$, we have thus contradicted Conjecture \ref{erdos} if $\eps$ is sufficiently small 
depending on $c_4$.  By modifying the above argument in a routine manner one thus obtains 
Conjecture \ref{furst-disc} for $c_3$ sufficiently small depending on $c_4$.

\section{The discretized Furstenburg Conjecture \ref{furst-disc} implies 
the Bilinear Distance Conjecture \ref{less-naive}}\label{furst-falcon}

To close the circle of implications and finish the proof of the Main Theorem 
\ref{grand-equiv} 
we need to show that the Discretized Furstenburg Conjecture \ref{furst-disc} implies 
the Bilinear Ring Conjecture \ref{less-naive}.  
This will be done by modifying the argument in Chung, Szemer\'edi, and Trotter 
\cite{cst:erdos}, in 
which the Szemer\'edi-Trotter theorem was applied to the discrete distance problem.  
The key geometric fact we use to pass from distances to lines is that if $|x_0 - x_1| = 
|x_0 - x_2|$, then $x_0$ lies on the perpendicular bisector of $x_1$ and $x_2$.  These 
lines need not point in different directions, but this will be remedied by a generic 
projective transformation.

\begin{figure}[htbp] \centering
\ \psfig{figure=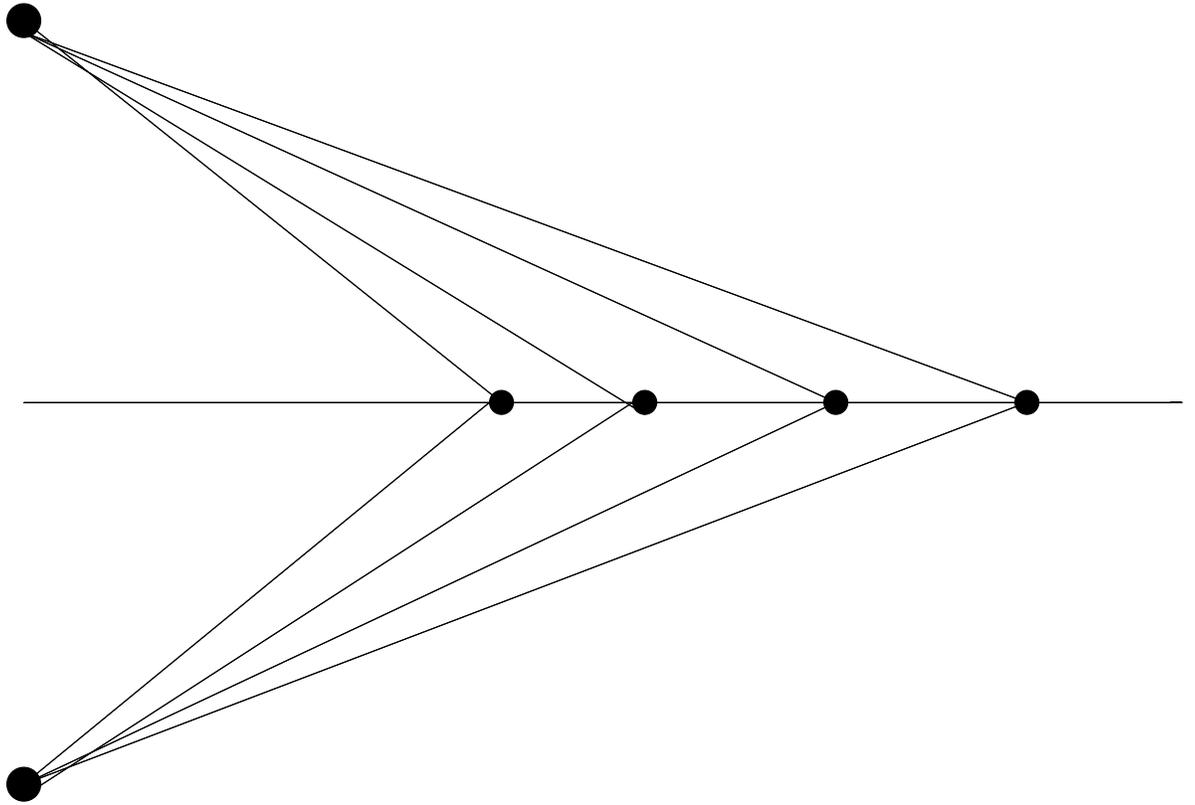}
\caption{Sets with few distances must concentrate on lines.}
\end{figure}

Assume that the Discretized Furstenburg Conjecture \ref{furst-disc} is true for some 
absolute constant $c_3 > 0$. Let 
$0 <  \eps \ll 1$ be fixed.  We may assume that $\delta$ is sufficiently small depending on 
$\eps$, since \eqref{eps1} is trivial otherwise.   Let $Q_j$, $E_j$, $D$ be as in the 
Bilinear Distance
Conjecture \ref{less-naive}.  Assume for contradiction that 
$$
|\{ (x_0,x_1,x_2) \in E_0 \times E_1 \times E_2: |x_0-x_1|, |x_0-x_2| \in D \}| \gtrapprox 
\delta^3
$$
We will obtain a contradiction from this, and it will be clear from the nature of the 
argument that \eqref{eps1} in fact holds for some absolute constant $c_1 > 0$ depending on 
$c_3$.

Let $E'_0$ denote the set of all $x_0 \in E_0$ such that
$$ | \{ x_2 \in E_2: |x_0 - x_2| \in D \} | \geq \delta^{\Ci \eps} \delta,$$
where $\Ci$ is an absolute constant to be chosen later.  We have
\be{dist-hyp} |\{ (x_0,x_1,x_2) \in E'_0 \times E_1 \times E_2: |x_0-x_1|, |x_0-x_2| 
\in D \}| \gtrapprox \delta^3
\end{equation}
if $\Ci$ is chosen to be sufficiently large and $\delta$ sufficiently small depending on 
$\Ci$ and $\eps$ (cf. \eqref{hyp3}.  Since the left hand side is clearly bounded by 
$|E'_0| |E_1| |E_2| \approx \delta^2 |E'_0|$, we thus see that $|E'_0| \approx \delta$.

Fix $\Ci$; all implicit constants may now depend on $\Ci$.  Since we clearly have
$$ | \{ x_2 \in E_2: |x_0 - x_2| \in D \} | \leq |E_2| \lessapprox \delta$$
then we see from \eqref{dist-hyp} that
$$ | \{ (x_0,x_1) \in E'_0 \times E_1: |x_0 - x_1| \in D \}| \gtrapprox \delta^2.$$
Since $D$ is $\delta$-discretized, we thus have
$$ | \{ (x_0,x_1,d) \in E'_0 \times E_1 \times D: |x_0 - x_1| \in D \cap \B(d,\delta) \}| 
\gtrapprox \delta^3.$$
Applying Cauchy-Schwarz \ref{cauchy} with $A = E_1 \times D$, $B = E'_0$, and $\lambda \approx 
\delta^{1/2}$ we thus have
$$ | \{ (x_0,x'_0,x_1,d) \in E'_0 \times E'_0 \times E_1 \times D: |x_0 - x_1|, 
|x'_0 - x_1| \in D \cap \B(d,\delta) \}| \gtrapprox \delta^{9/2}.$$
For fixed $x_0, x'_0, x_1$, the set of all $d$ which contribute to the above set has 
measure $O(\delta)$, and vanishes unless $|x_0 - x_1| = |x'_0 - x_1| + O(\delta)$.  
Thus we have
$$ | \{ (x_0,x'_0,x_1) \in E'_0 \times E'_0 \times E_1: |x_0-x_1|, |x'_0 - x_1| \in D, 
|x_0 - x_1| = |x'_0 - x_1| + O(\delta) \}| \gtrapprox \delta^{7/2}.$$

Let $\Cj$ be an absolute constant to be chosen later.

\begin{lemma} We have
\bas  | \{ &(x_0,x'_0,x_1) \in E'_0 \times E'_0 \times E_1: |x_0 - x_1|, |x'_0 - x_1| 
\in D, |x_0 - x_1| = |x'_0 - x_1| + O(\delta);\\
& |x_0 - x'_0| \leq \delta^{\Cj \eps} \}| \lessapprox \delta^{\Cj \eps/2} \delta^{7/2}.
\end{align*}
\end{lemma}

\begin{proof}  Since $|E'_0|, |E_1| \approx \delta$, it suffices to show that
$$  | \{ x_0 \in E'_0: |x_0 - x_1| = |x'_0 - x_1| + O(\delta); |x_0 - x'_0| \leq 
\delta^{\Cj \eps} \}| \lessapprox \delta^{\Cj \eps/100} \delta^{3/2}$$
for all $x'_0$, $x_1$ in $E'_0$, $E_1$ respectively.

Fix $x'_0$, $x_1$.  By definition of $E'_0$ it suffices to show that
$$  | \{ (x_0, x_2) \in E'_0 \times E_2: |x_0 - x_1| = |x'_0 - x_1| + O(\delta); 
|x_0 - x'_0| \leq \delta^{\Cj \eps}; |x_0 - x_2| \in D \}| \lessapprox 
\delta^{\Cj \eps/100} \delta^{5/2}.$$
Since $|E_2| \approx \delta$, it thus suffices to show that
\be{tar}  | \{ x_0 \in E'_0: |x_0 - x_1| = |x'_0 - x_1| + O(\delta); |x_0 - x'_0| \leq 
\delta^{\Cj \eps}; |x_0 - x_2| \in D \}| \lessapprox
\delta^{\Cj \eps/100} \delta^{3/2}
\end{equation}
for all $x_2 \in E_2$.

Fix $x_2$.  The set in \eqref{tar} is contained in an annular arc of thickness $O(\delta)$, 
angular width $O(\delta^{\Cj \eps)}$, and radius $\approx 1$ centered at $x_1$.  From 
\eqref{gen-sep} and elementary geometry that the possible values of $|x_0 - x_2|$ thus 
lie in an interval of length $\lessapprox \delta^{\Cj \eps}$.  Since $D$ is a 
$(\delta,\frac{1}{2})_1$ set, we thus see that the possible values of $|x_0 - x_2|$ are 
contained in the union of $\lessapprox \delta^{\Cj \eps/2} \delta^{-1/2}$ intervals of 
length $\delta$.  From \eqref{gen-sep} and elementary geometry we thus see that the set 
in \eqref{tar} can be covered by $\lessapprox \delta^{\Cj \eps/2} \delta^{-1/2}$ balls of 
radius $\delta$.  The claim follows.
\end{proof}

If we choose $\Cj$ to be sufficiently large, and $\delta$ is sufficiently small depending on 
$\Cj$ and $\eps$, we thus have from the above that
$$  | \{ (x_0,x'_0,x_1) \in E'_0 \times E'_0 \times E_1: |x_0 - x_1|, |x'_0 - x_1| \in D; 
|x_0 - x_1| = |x'_0 - x_1| + O(\delta); |x_0 - x'_0| \approx 1 \}| \gtrapprox  
\delta^{7/2},$$
where the implicit constants can now depend on $\Cj$.  Since $|E'_0| \approx \delta$, we 
can thus find $x_0 \in E'_0$ such that
$$  | \{ (x'_0,x_1) \in E'_0 \times E_1: 
|x_0 - x_1|, |x'_0 - x_1| \in D; |x_0 - x_1| = |x'_0 - x_1| + O(\delta); |x_0 - x'_0| 
\approx 1 \}| \gtrapprox \delta^{5/2}.$$

Fix this $x_0$.  Let $E''_0$ denote the set
$$ E''_0 := \{ x'_0 \in E'_0: |x_0 - x'_0| \approx 1 \}.$$
For each $x'_0 \in E''_0$, let $R[x'_0]$ denote the set
$$ R[x'_0] := \{ x_1 \in E_1: |x_0 - x_1| \in D; |x_0 - x_1| = |x'_0 - x_1| + O(\delta) \}.$$
We thus have
\be{int}
 \int_{E''_0} |R[x'_0]|\ dx'_0 \gtrapprox \delta^{5/2}.
\end{equation}
From elementary geometry, and the fact that $Q_0$ and $Q_1$ have separation $\approx 1$, we 
see that $R[x'_0]$ is contained in a rectangle $\overline{R}[x'_0]$ of dimensions 
$\approx 1 \times \delta$.  Since $D$ is a $(\delta,\frac{1}{2})_1$ set, it is easy to see 
from elementary geometry that $R[x'_0]$ is a $(\delta, \frac{1}{2})_2$ set.  In particular, 
we have $|R[x'_0]| \lessapprox \delta^{3/2}$, which implies from \eqref{int} that 
$|E''_0| \gtrapprox \delta$.
Since $|E_0| \approx \delta$, we thus have $|E''_0| \approx \delta$.  From 
\eqref{int} 
again we thus see that
$$ |\{ x'_0 \in E''_0: |R[x'_0]| \approx \delta^{3/2} \}| \approx \delta$$
for suitable choices of constants.  In particular, we can find a $\delta$-separated set 
$\Sigma \subset E''_0$ such that $\# \Sigma \approx \delta^{-1}$ and $|R[x'_0]| \approx 
\delta^{3/2}$ for all $x'_0 \in \Sigma$.

The sets $R[x'_0]$ resemble the sets $R_\omega$ in the hypothesis of Conjecture 
\ref{furst-disc}, but their orientations need not be $\delta$-separated.  To remedy this we 
apply a projective linear transformation.

To find the right transformation to use, we first must isolate a line in $\R^2$ in which 
the (extensions of) $\overline{R}[x'_0]$ are well-separated.  To make this precise we apply 
some normalizations.  By a rescaling we may let $Q_1$ be the square $[0,1] \times [0,1]$, 
and by a refinement we may assume that the the direction of the $\overline{R}[x'_0]$ are 
within $\pi/4$ of the direction $(1,0)$.  In particular, if we extend the long side of 
$\overline{R}[x'_0]$ to have length $10$, it will intersect the strip $[2,3] \times \R$ in 
a parallelogram $P[x'_0]$ of thickness $\approx \delta$ and slope $O(1)$.

We now apply 

\begin{lemma}\label{pseudo-cordoba}  We have
$$ \| \sum_{x'_0 \in \Sigma} \chi_{P[x'_0]} \|_2^2 \lessapprox 1.$$
\end{lemma}

\begin{proof} We shall use C\'ordoba's argument, using \eqref{dim} for $E_0$ as a substitute 
for the direction-separation property.  Expand out the left-hand side as
$$ \sum_{x'_0 \in \Sigma} \sum_{x''_0 \in \Sigma} |P[x'_0] \cap P[x''_0]|.$$
Since $\# \Sigma \approx \delta^{-1}$, it thus suffices to show that
$$ \sum_{x'_0 \in \Sigma} |P[x'_0] \cap P[x''_0]| \lessapprox \delta$$
for all $x''_0 \in \Sigma$.

Fix $x''_0$.  The quantity $|P[x'_0] \cap P[x''_0]|$ can vary from 0 to $\approx \delta$.  
We need only consider the contribution when $\delta^2 \lessapprox |P[x'_0] \cap P[x''_0]| 
\lessapprox \delta$, since the remaining contribution is trivial to handle.  By dyadic 
pigeonholing and absorbing the logarithmic factor into the $\lessapprox$ symbol, it suffices 
to show the distributional estimate
$$ \# \{x'_0 \in \Sigma: |P[x'_0] \cap P[x''_0]| \approx \sigma \delta \} \lessapprox 
\frac{1}{\sigma}$$
for all dyadic $\delta \lessapprox \sigma \lessapprox 1$.

Fix $\sigma$.  The set $P[x'_0]$ lies within $\lessapprox \delta$ of the perpendicular 
bisector of $x'_0$ and $x_0$, and within $\lessapprox 1$ of $x'_0$ and $x_0$, which are 
themselves separated by $\approx 1$.  Similarly for $P[x''_0]$.  From elementary geometry 
we thus see that
$$ |P[x'_0] \cap P[x''_0]| \approx \sigma \delta \implies |x'_0 - x''_0| \lessapprox 
\delta/\sigma.$$
Since $x'_0, x''_0$ lie within a $\delta$-separated subset of $E_0$, which is a 
$(\delta,1)_2$ set, we see from \eqref{dim} that
$$ \# \{ x'_0 \in \Sigma: |x'_0 - x''_0| \lessapprox \delta/\sigma \} \lesssim  
\frac{1}{\sigma}.$$
The claim follows.
\end{proof}

To complement this $L^2$ bound we have the trivial $L^1$ bound
$$ \| \sum_{x'_0 \in \Sigma} \chi_{P[x'_0]} \|_1 = \sum_{x'_0 \in \Sigma} |P[x'_0]| 
\approx \delta^{-1} \delta = 1.$$
From H\"older's inequality we thus see that $\sum_{x'_0 \in \Sigma} \chi_{P[x'_0]}$ must be 
supported on a set of measure $\gtrapprox 1$, so that
$$ |\bigcup_{x'_0 \in \Sigma} P[x'_0]| \gtrapprox 1.$$
Since the set in the left-hand side is contained in the strip $[2,3] \times R$, we can 
thus find a $2 \leq \x \leq 3$ such that
$$ |\{ \y: (\x,\y) \in \bigcup_{x'_0 \in \Sigma} P[x'_0] \}| \gtrapprox 1.$$
Fix this $\x$.  Each $x'_0 \in \Sigma$ contributes an interval of length $\approx \delta$ to 
the above set.  Thus we can find a refinement $\Sigma'$ of $\Sigma$ such that the sets
$\{ \y: (\x,\y) \in P[x'_0] \}$ 
are separated by $\gg \delta$.

Let $L$ be a projective transformation which sends the line $\{\x\} \times \R$ to the line 
at infinity, but maps $[0,1] \times [0,1]$ to a bounded set and has Jacobian $\approx 1$ on 
$[0,1] \times [0,1]$.  Thus $L(E_1)$ is a $(\delta,1)_2$ set with measure $\approx \delta$.

For each $x'_0 \in \Sigma'$, we see from elementary projective geometry we see that the 
sets $L(R[x'_0])$ are $(\delta,\frac{1}{2})_2$ sets contained in a rectangle of dimensions 
$\approx 1 \times \delta$, and the orientation $\omega = \omega(x'_0)$ of these rectangles 
are $\delta$-separated as $x'_0$ varies along $\Sigma'$.  Write $\Omega$ for the set of all 
the orientations $\omega$ arising in this manner, so that $\# \Omega \approx \delta^{-1}$, 
and write $R_\omega := L(R[x'_0])$ for all $x'_0 \in \Sigma'$.  Also write $E := L(E_1)$.  
Since $R_\omega$ is a $(\delta,\frac{1}{2})_2$ set and $|R_\omega| \approx \delta^{3/2}$,
 we have
$$ |\{ (x_0,x_1) \in R_\omega: |x_0 - x_1| \approx 1 \}| \gtrapprox \delta^3$$
for appropriate choices of constants.  For any fixed $x_0, x_1$ with $|x_0 - x_1| 
\approx 1$, there are at most $\lessapprox 1$ values of $\omega$ for which $(x_0,x_1)$ is 
contained in the above set, thanks to the $\delta$-separation of the $\omega$.  Since 
$R_\omega \subset E$, we may thus sum the above estimate in $\omega$ to obtain
$$ |\{ (x_0,x_1) \in E: x_0, x_1 \in R_\omega \hbox{ for some } \omega \in \Omega; 
|x_0 - x_1| \approx 1 \}| \gtrapprox \delta^2.$$
But this contradicts \eqref{furst-est} if $\eps$ is sufficiently small depending on $c_3$ 
and $\delta$ is sufficiently small depending on $\eps$.  By modifying 
the above argument in a routine manner one thus obtains the
Bilinear Distance Conjecture \ref{less-naive} for 
$c_1$ sufficiently small depending on $c_3$.  This completes the proof of the Main Theorem 
\ref{grand-equiv}.

\section{Discretization of fractals}\label{discretization}

In order to pass from the $\delta$-discretized Bilinear Distance
conjecture \ref{less-naive} and Discretized Furstenburg conjecture
\ref{furst-disc} to their respective continuous analogues the Distance
Conjecture \ref{dist-conj} and the Furstenburg conjecture
\ref{furst-conj} we will need some tools to cover an $\alpha$-dimensional set in 
$\R^n$ by $(\delta,\alpha+C\eps)_n$ sets for various values of $\delta$.

We begin this section by recalling the definition of Hausdorff dimension.

\begin{definition}  Let $\alpha > 0$.  For any bounded set $E$ and $c > 0$, we define the 
\emph{Hausdorff content} $h_{\alpha,c}(E)$ to be the infimum of the quantity
$$ \sum_{i \in I} r_i^\alpha$$
where $\{ \B(x_i,r_i) \}_{i \in I}$ ranges over all collections of balls of radii 
$r_i < c$ which cover $E$. 
$$ \dim(E) := \inf \{ \alpha: \inf_{c>0} h_{\alpha,c}(E) = 0 \} = \sup \{ \alpha: 
\sup_{c>0} h_{\alpha,c}(E) = +\infty \}.$$
\end{definition}

\begin{definition}  Let $\{ X_\alpha \}_{\alpha \in A}$ be a countable collection of sets.  
We say that the $X_\alpha$ \emph{strongly cover} $E$ if each point in $E$ is contained in 
infinitely many sets $X_\alpha$.
\end{definition}

We shall require a variant of the Borel-Cantelli lemma for Hausdorff content. 

\begin{lemma}\label{bc}  Let $0 < \alpha \leq n$, and let $X_i \subset \R^n$ for $i \in 
\Z$ be such that
$$\sum_{i=1}^\infty h_{\alpha,c}(X_i) < \infty$$
for some $c > 0$.  Suppose also that the $X_i$ strongly cover a set $E$.  Then $\dim(E) 
\leq \alpha$.
\end{lemma}

\begin{proof}  
For any integer $N$, we have $E \subset \bigcup_{i > N} X_i$.  Since Hausdorff content is 
sub-additive, we thus have
$$ h_{\alpha,c}(E) \leq \sum_{i > N} h_{\alpha,c}(X_i).$$
The claim then follows by letting $N \to \infty$.
\end{proof}

We can now prove a covering lemma, which is the main result of this section.
For technical reasons it will be convenient to not work with dyadic $\delta$ as we have done 
in the past, but move to a much sparser range of scales, namely the hyper-dyadic scales 
(cf. \cite{borg:high-dim}).  More precisely:

\begin{definition}  Let $0 < \eps \ll 1$ be given.  We call a number \emph{hyper-dyadic} if 
it is of the form $2^{-\floor{(1+\eps)^k}}$ for some integer $k \geq 0$, where $\floor{x}$ 
is the integer part of $x$.  We call a cube \emph{hyper-dyadic} if it is dyadic and 
its side-length is hyper-dyadic.
\end{definition}

Note that there are at most $C_\eps$ hyper-dyadic numbers between $\delta$ and 
$\delta^{100}$ for any choice of $\delta$, in contrast to $C \log(1/\delta)$ in the dyadic 
regime.  This improved bound will be important in the proof of Theorem \ref{falc-reduc}.

\begin{lemma}\label{covering}  Let $0 < \eps \ll 1$, $0 < \alpha < n$, and let $E$ be a 
compact subset of $\R^n$.
\begin{itemize}
\item If $\dim(E) \leq \alpha$, then one can associate a $(\delta,\alpha)_n$ set $X_\delta$ 
to each hyper-dyadic $\delta$ such that the $X_\delta$ strongly cover $E$.
\item Conversely, if $C$ is sufficiently large and there is a $(\delta,\alpha-C\eps)_n$ set 
$X_\delta$ for each hyper-dyadic $\delta$ such that the $X_\delta$ strongly cover $E$, then 
$\dim(E) \leq \alpha$.
\end{itemize}
\end{lemma}

\begin{proof}  We first prove the latter claim.  Since $X_\delta$ is a
$(\delta,\alpha-C\eps)_n$ set we can cover it (if the constants are chosen appropriately) 
by about $\delta^{-\alpha+\eps}$ balls of radius $\delta$, so that
$$ h_{\alpha,1}(X_\delta) \leq C \delta^\eps.$$
The claim then follows from Lemma \ref{bc}.

Now we show the former claim.
Fix $E$.  For every hyper-dyadic number $c$, we can find a collection 
$\{ \B(x_{c,i},r_{c,i}) \}_{i \in I_c}$ of balls covering $E$ such that 
$r_{c,i} < c$ and
\be{ycr} \sum_{i \in I_c} r_{c,i}^{\alpha+C\eps} \ll 1.
\end{equation}
By reducing the constant $C$ slightly we may assume that the $r_{c,i}$ are hyper-dyadic.  

For each hyper-dyadic $r$, let $Y_{c,r}$ denote the set
$$ Y_{c,r} := \bigcup_{i \in I_c: r_{c,i} = r} \B(x_{c,i},r_{c,i}).$$
Clearly the sets $Y_{c,r}$ strongly cover $E$ as $c$, $r$ both vary.

Fix $c$, $r$, and let $\Q_{c,r}$ be a collection of hyper-dyadic cubes $Q$ of side-length 
at least $r$ which cover $Y_{c,r}$ and which minimize the quantity
$$ \sum_{Q \in \Q_{c,r}} l(Q)^\alpha,$$
where $l(Q)$ denotes the side-length of $Q$.
Such a minimizer exists since there are only a finite number of hyper-dyadic cubes which 
are candidates for inclusion in $\Q_{c,r}$.  From \eqref{ycr} one can cover $Y_{c,r}$ by at 
most $r^{-\alpha-\eps}$ cubes of side-length $r$, hence
\be{lq-card} \sum_{Q \in \Q_{c,r}} l(Q)^\alpha \leq C r^{-\eps}.
\end{equation}
In particular, we have $l(Q) \leq C r^{-\eps/\alpha}$ for all $Q \in \Q_{c,r}$.

From the construction of $\Q_{c,r}$ we see that the $Q$ are all disjoint, and
for all hyper-dyadic cubes $I$ we have
\be{dyadic-dim}
 \sum_{Q \in \Q_{c,r}: Q \subset I} l(Q)^\alpha \leq l(I)^\alpha
\end{equation}
since we could otherwise remove those cubes in $I$ from $\Q_{c,r}$ and replace them with 
$I$, contradicting minimality.

For each dyadic $r \leq \delta \leq C r^{-\eps/\alpha}$, let $X_{\delta,c,r}$ denote the
 set
$$ X_{\delta,c,r} := (\bigcup_{Q \in \Q_{c,r}: l(Q) = \delta} Q) + \B(0,\delta).$$
Clearly $X_{\delta,c,r}$ is a $\delta$-discretized set.  From \eqref{dyadic-dim} we see that 
$X_{\delta,c,r}$ is in fact a $(\delta,\alpha)_n$ set.

Now define $X_\delta := \bigcup_{c,r} X_{\delta,c,r}$.  From the constraints
$r < c$ and $\delta < C r^{-\eps/\alpha}$ we see that there are at most $C \log(1/\delta)^2$ 
pairs $(c,r)$ associated to each $\delta$.  Hence $X_\delta$ is also a $(\delta,\alpha)_n$ 
set.  By construction we see that the $X_\delta$ strongly cover $E$, and so we are done.
\end{proof}

\section{The Discretized Furstenburg conjecture implies the Furstenburg
problem.}\label{furst-disc-sec}

We now prove Theorem \ref{furst-equiv}.  Suppose that the
Discretized Furstenburg Conjecture \ref{furst-disc} holds for 
some $c_3 > 0$, and let $K$ be a $\frac{1}{2}$-set using the notation of the introduction.  
Let $0 < \eps \ll c_2 \ll c_3^2$ be constants to be chosen later.  Assume for contradiction 
that $K$ has Hausdorff dimension less than $1 + c_2$.

By Lemma \ref{covering}, we may find a $(\delta,1 + c_2)_2$ set $X_\delta$ for each 
hyper-dyadic $\delta$ such that the $X_\delta$ strongly cover $K$.

If $\omega \in S^1$ is a direction, we call $\omega$ \emph{bad with respect to $\delta$} if 
one can find a line $l$ in the direction $\omega$ such that 
\be{hbad}
h_{1/2 - c_2,1}(l \cap X_\delta) \geq \delta^{c_2},
\end{equation}

The main estimate we need is 

\begin{lemma}  For all hyper-dyadic $\delta$, we have
\be{key-furst}
 |\{ \omega \in S^1: \omega \hbox{ is bad with respect to } \delta \}| \lessapprox C_{c_2} 
\delta^{Cc_2}
\end{equation}
if $c_2$ is sufficiently small with respect to $c_3^2$.
\end{lemma}

\begin{proof}
The proof of trivial if $\delta$ is large, so we will assume that $\delta$ is sufficiently 
small depending on $c_2$, $c_3$.

From Kakeya \ref{cordoba} with $f := \chi_{X_\delta}$ and Chebyshev's inequality we have
$$ \{ \omega \in S^1: (\chi_{X_\delta})^*_\delta(\omega) > \delta^{-c_2} \delta^{1/2} \} 
\lessapprox \delta^{Cc_2}.$$
Thus to show \eqref{key-furst} it suffices to show that
\be{key-furst2}
|\{ \omega \in S^1: \omega \hbox{ is bad with respect to } \delta, 
(\chi_{X_\delta})^*_\delta(\omega) < \delta^{-c_2} \delta^{1/2} \}| \lessapprox C_{c_2} 
\delta^{Cc_2}.
\end{equation}

Suppose for contradiction that \eqref{key-furst2} failed.  Let $\Omega$ be a maximal 
$\delta$-separated subset in the set in \eqref{key-furst2}; we thus have
\be{car}
\# \Omega \gtrapprox C_{c_2} \delta^{Cc_2} \delta^{-1}.
\end{equation}
By construction, for each $\omega \in \Omega$ we can find a line $l_\omega$ in the direction 
$\omega$ such that \eqref{hbad} holds.  Let $R_\omega$ denote the set 
$(l_\omega + \B(0,\delta)) \cap X_\delta$.  From the construction of $\Omega$ we thus have
\be{rom-bound}
 |R_\omega| \lessapprox \delta (\chi_{X_\delta})^*_\delta(\omega)
\lessapprox \delta^{-Cc_2} \delta^{3/2}.
\end{equation}

Let $\Q$ be a collection of squares $Q$ of side-length $l(Q) \geq \delta$ which covers 
$R_\omega$ and which minimizes the quantity
$$ \sum_{Q \in \Q} l(Q)^{1/2 - \sqrt{c_2}}.$$
As in the proof of Lemma \ref{covering}, a minimizer $\Q$ exists and the squares in 
$Q$ are disjoint and satisfy
\be{dyadic-dim2}
 \sum_{Q \in \Q: Q \subset I} l(Q)^{1/2 - \sqrt{c_2}} \leq l(I)^{1/2 - \sqrt{c_2}}
\end{equation}
for all squares $I$.
Also, for all $Q \in \Q$ we have
$$ |Q \cap R_\omega| \gtrapprox \delta^2 (l(Q)/\delta)^{1/2 - \sqrt{c_2}}$$
since otherwise we could replace $Q$ by all the $\delta$-cubes contained in $Q$, 
contradicting the minimality of $\Q$.  Summing this over all $Q$ we obtain
$$ |R_\omega| \gtrapprox \delta^{3/2 + \sqrt{c_2}} \sum_{Q \in \Q} l(Q)^{1/2 - 
\sqrt{c_2}}.$$
From \eqref{rom-bound} we thus obtain
$$ \sum_{Q \in \Q} l(Q)^{1/2 - \sqrt{c_2}} \lessapprox \delta^{-\sqrt{c_2}-C c_2}.$$
We thus have
$$ \sum_{Q \in \Q: l(Q) > \delta^{1 - A\sqrt{c_2}} } l(Q)^{1/2 - c_2} \lessapprox 
\delta^{(1-A\sqrt{c_2})(\sqrt{c_2}-c_2)} \delta^{-\sqrt{c_2} - C c_2}$$
If we choose $A$ sufficiently large, we thus have (for $\delta$ sufficiently small)
$$ \sum_{Q \in \Q: l(Q) > \delta^{1 - A\sqrt{c_2}} } l(Q)^{1/2 - c_2} \ll \delta^{c_2}.$$
In particular, we have
$$ h_{1/2 - c_2,1}( \bigcup_{Q \in \Q: l(Q) > \delta^{1 - A \sqrt{c_2}}} Q )
\ll \delta^{c_2}.$$
On the other hand, since $\omega$ is bad with respect to $\delta$, we have
$$ h_{1/2 - c_2,1}( R_\omega ) \geq h_{1/2 - c_2,1}(l_\omega \cap X_\delta) 
\geq \delta^{c_2}.$$
Thus, if we let $R'_\omega$ denote the set
$$ R'_\omega := (R_\omega \backslash \bigcup_{Q \in \Q: l(Q) > \delta^{1 - A \sqrt{c_2}}} 
Q) + \B(0,\delta)$$
then we have
$$ h_{1/2 - c_2,1}( R'_\omega ) \gtrapprox \delta^{c_2}.$$
Since $R'_\omega$ is $\delta$-discretized, we have in particular that
$$ |R'_\omega| \gtrapprox \delta^{3/2 + C c_2}.$$
The set $R'_\omega$ is covered by the dilates of those cubes $Q \in \Q$ for which $l(Q) 
\leq \delta^{1 - A \sqrt{c_2}}$.  From this and \eqref{dyadic-dim2} we see that 
$R'_\omega$ is a $(\delta,1/2)_2$ set but with $\eps$ replaced by $A \sqrt{c_2}$.  
From \eqref{furst-est} we thus have
$$
| \{ (x_0,x_1) \in K \times K: x_1, x_0 \in R'_{\omega} \hbox{ for some } 
\omega \in \Omega \} | \lessapprox \delta^{2+c_3-C \sqrt{c_2}}.
$$
On the other hand, from Separation \ref{sep-lemma} we have
$$ | \{ (x_0, x_1) \in R'_\omega: |x_0 - x_1| \gtrapprox \delta^{C A \sqrt{c_2}} \} | 
\gtrapprox \delta^{3 + C A \sqrt{c_2}}.$$
Summing this on $\omega$ using \eqref{car} and noting that each $(x_0,x_1)$ can be in at 
most $\lessapprox \delta^{-C A \sqrt{c_2}}$ of the above sets, we obtain
$$ | \{ (x_0,x_1) \in K \times K: x_1, x_0 \in R'_{\omega} \hbox{ for some } 
\omega \in \Omega \} | \gtrapprox C_{c_2} \delta^{2+ C A \sqrt{c_2}}.$$
If $c_2$ is sufficiently small with respect to $c_3^2$ we obtain the desired contradiction, 
if $\delta$ is sufficiently small.
\end{proof}

If $\eps$ is chosen sufficiently small depending on $c_2$, then the left-hand side of 
\eqref{key-furst} is thus summable in $\delta$.  By the Borel-Cantelli lemma (for Lebesgue 
measure) we can thus find a direction $\omega$ which is only bad with respect to a finite 
number of hyper-dyadic $\delta$.  In particular we have
$$
\sum_\delta h_{1/2 - c_2,1}(l \cap X_\delta) < \infty$$
for all lines $l$ parallel to $\omega$.  Since the $l \cap X_\delta$ strongly covers 
$l \cap K$, we thus see from Lemma \ref{bc} that $\dim(l \cap K) < 1/2$ for all $l$ parallel 
to $\omega$.  But this contradicts the assumption that $K$ is a $\frac{1}{2}$-set.
\endprf

We remark that a similar result obtains for all $\beta$-sets providing that $\beta$ is 
sufficiently close to $\frac{1}{2}$ (depending on $c_3$).

\section{The Bilinear Distance conjecture implies the Falconer
Distance Conjecture}\label{falc-reduc-sec}

We now prove Theorem \ref{falc-reduc}.  Suppose that the Bilinear
Distance Conjecture \ref{less-naive} holds for 
some $c_1 > 0$.  Let $0 < \eps \ll c_0 \ll c_1$ be constants to be chosen later.  Assume for 
contradiction that one can find a compact set $K$ with dimension $\dim(K) \geq 1$ such 
that $\dim (\dist(K)) \leq 1/2 + c_0$.

By Frostman's lemma \cite{mat:GMT} we may find a probability measure $\mu$ supported on $K$ such 
that 
\be{frosty}
\mu(\B(x,r)) \leq C_{\eps} r^{1-\eps}
\end{equation}
for all balls $\B(x,r)$.  Fix this $\mu$.

By Lemma \ref{covering} we may find a $(\delta,1/2 + c_0)_1$ set $D_\delta$ for each 
hyper-dyadic $\delta$ such that the $D_\delta$ strongly cover $\dist(K)$.  If one then 
defines
$$ X_\delta := \{ (x,y) \in K \times K: |x-y| \in D_\delta \}$$
then the $X_\delta$ strongly cover $K \times K$.

If it were not for the bilinear formulation of Conjecture \ref{less-naive}, one could hope 
to prove a bound like 
\be{mux}
\mu(X_\delta) \lessapprox \delta^{C^{-1} c_0},
\end{equation}
if $c_0 \ll c_1$, which would allow us to obtain a contradiction from the Borel-Cantelli 
lemma.  These types of bounds however are not achievable because of the counter-example 
\eqref{counter}.  Furthermore, it is possible for $K$ to contain obstructions like 
\eqref{counter} at infinitely many scales.  Fortunately, one can show that it is not 
possible for too many pairs $x,y \in K$ to be simultaneously contained in sets like 
\eqref{counter} at infinitely many scales, which allows one to proceed.  Of course, one 
has to make the notion of ``looking like \eqref{counter}'' precise, which causes some 
unpleasant technicalities.

We begin by converting \eqref{eps1} to a bilinear variant of \eqref{mux}.

\begin{lemma}\label{est-useful}  If $c_0$ is sufficiently small depending on $c_1$, and 
$\delta$ is sufficiently small, then
\be{bilinear-est}
\mu^3(\{ (x_0,x_1,x_2): (x_0,x_1), (x_0,x_2) \in X_\delta; 
|(x_1 - x_0) \wedge (x_2 - x_0)| \geq \delta^{C^{-1} c_0} \})
\leq C_{\eps,c_0} \delta^{C^{-1} c_0}.
\end{equation}
where $\mu^3$ is product measure on $K \times K \times K$.
\end{lemma}

\begin{proof}
Let $c_5$ be a constant to be chosen later.
We shall show that
\be{bil-precise}
\begin{split}
\mu^3(\{ (x_0,x_1,x_2): &(x_0,x_1), (x_0,x_2) \in X_\delta; 
|(x_1 - x_0) \wedge (x_2 - x_0)| \geq \delta^{c_5} \})\\
&\leq C_{\eps,c_0} (\delta^{C^{-1} c_5} + \delta^{-C c_5} \delta^{C^{-1} c_0}),
\end{split}
\end{equation}
from which \eqref{bilinear-est} follows from a suitable choice of $c_5$.

Partition $K = K_1 \cup K_2$, where
$$ K_1 := \{ x \in K: \mu(B(x,\delta)) \geq \delta^{c_5} \delta \},$$
$$ K_2 := \{ x \in K: \mu(B(x,\delta)) < \delta^{c_5} \delta \}.$$
Let us first deal with the contribution to \eqref{bil-precise} of the case when at least one 
of $x_0, x_1, x_2$ is in $K_2$.  By Fubini's theorem and symmetry it suffices to show that
$$
\mu^2(\{ (x_0,x_1): x_0 \in E_2, (x_0,x_1) \in X_\delta \})
\leq C_{\eps,c_0} \delta^{c_5/10}.
$$

By Cauchy-Schwarz \ref{cauchy} it suffices to show that
$$
\mu^3(\{ (x_0,x_1,x_2): x_0 \in E_2, (x_0,x_1), (x_0,x_2) \in X_\delta \})
\lessapprox \delta^{c_5/5}.
$$

Let us first consider the contribution of the case $|x_1 - x_2| \lessapprox \delta^{c_5/5}$.  For each $x_0$,
 $x_1$, the set of $x_2$ which contribute to the above expression has measure $O(\delta^{c_5/5})$ by \eqref{frosty}, and so this contribution is acceptable by Fubini's theorem.  It thus remains to show
$$
\mu^3(\{ (x_0,x_1,x_2): x_0 \in E_2, (x_0,x_1), (x_0,x_2) \in X_\delta, |x_1 - x_2| 
\gtrapprox \delta^{c_5/5} \})
\lessapprox \delta^{c_5/5}.
$$
By Fubini's theorem again, it suffices to show that
$$
\mu(\{ x_0 \in E_2: (x_0,x_1), (x_0,x_2) \in X_\delta\})
\lessapprox \delta^{c_5/5}
$$
for all $x_1, x_2 \in E$ such that $|x_1 - x_2| \gtrapprox \delta^{c_5/5}$.

Fix $x_1$, $x_2$.  By \eqref{frosty} again, it suffices to show that
$$
\mu(\{ x_0 \in E_2: |x_0 - x_1|, |x_0 - x_2| \in D_\delta; |x_0 - x_1|, |x_0 - x_2| 
\gtrapprox \delta^{c_5/5} \})
\lessapprox \delta^{c_5/5}.
$$
The set $D_\delta$ can be covered by $\lessapprox \delta^{-1/2}$ intervals of length 
$\delta$.  Let $\I$ denote the collection of those intervals $I$ in this cover such that 
$\dist(0,I) \gtrapprox \delta^{c_5/5}$.  It suffices to show that
\be{ij-sum}
\sum_{I,J \in \I} \mu(\{ x_0 \in E_2: |x_0 - x_1| \in I; |x_0 - x_2| \in J \})
\lessapprox \delta^{c_5/5}.
\end{equation}
For fixed $I$, $J$, the set described above is contained in an annular arc of thickness 
$\delta$, radius $\delta^{c_5/5} \lessapprox r \lessapprox 1$, and angular width bounded 
by
$$ \lessapprox \frac{\delta^{-c_5/10} \delta}{(\delta + | \dist(I,J) - |x_1-x_2| |)^{1/2}}$$
as one can easily compute using elementary geometry.  By the construction of $E_2$ and a
simple covering argument, we can thus bound the left-hand side of \eqref{ij-sum} by
$$
\lessapprox \delta^{c_5} \frac{\delta^{-c_5/10} \delta}{(\delta + 
| \dist(I,J) - |x_1-x_2| |)^{1/2}}.$$
To complete the proof of \eqref{ij-sum} it thus suffices to show that
$$ \sum_{J \in \I} \frac{\delta^{-c_5/10} \delta}{(\delta + | \dist(I,J) - 
|x_1-x_2| |)^{1/2}} \lessapprox \delta^{1/2}$$
for all $I \in \I$.  But this follows easily by dyadically decomposing the $J$ based on 
$\dist(I,J)$ and noting from \eqref{dim} that for each $k$, there are $\lessapprox 2^{k/2}$ 
intervals $J$ for which $\dist(I,J) \approx 2^k \delta$.
Note that any logarithmic factors can be absorbed into the $\lessapprox$ notation.

To conclude the proof of \eqref{bil-precise} it remains to show that
\bas
\mu^3(\{ &(x_0,x_1,x_2) \in K_1 \times K_1 \times K_1: (x_0,x_1), (x_0,x_2) \in X_\delta; 
\\
&|(x_1 - x_0) \wedge (x_2 - x_0)| \geq \delta^{c_5} \})
\leq C_{\eps,c_0} \delta^{-C c_5} \delta^{C^{-1} c_0}.
\end{align*}
From the definition of $K_1$ and \eqref{frosty} we see that $K_1$ is contained in a 
$(\delta,\frac{1}{2} + c_5)_2$ set $E$.  From \eqref{frosty} and a covering argument we 
have
\bas
\mu^3(\{ &(x_0,x_1,x_2) \in K_1 \times K_1 \times K_1: (x_0,x_1), (x_0,x_2) \in X_\delta; 
|(x_1 - x_0) \wedge (x_2 - x_0)| \geq \delta^{c_5} \})\\
&\lessapprox \delta^{-3}
|\{ (x_0,x_1,x_2) \in E \times E \times E: |x_0-x_1|, |x_0-x_2| \in D_\delta + 
B(0,\delta^{1-C\eps});\\
&|(x_1 - x_0) \wedge (x_2 - x_0)| \geq \frac{1}{2} \delta^{c_5} \}|.
\end{align*}
The claim then follows from \eqref{eps1}.
\end{proof} 

Henceforth we assume that $c_0$ is so small that Lemma \ref{est-useful} holds.

We shall use \eqref{bilinear-est} to create a dichotomy, that either \eqref{mux} holds or 
that the pairs in $X_\delta$ are concentrated in a thin set resembling \eqref{counter}.  
More precisely, we have

\begin{lemma}\label{excep}  
If $\Cl$ is a sufficiently large constant, then for each $\delta$ we can find an integer 
$N_\delta$ and sets $S_{\delta,1}, \ldots, S_{\delta,N_\delta} \subset B(0,C)$ such that
\begin{itemize}
\item For each $\delta$, each $x \in K$ is contained in at most $C$ sets $S_{\delta,i}$, 
where $C$ is an absolute constant independent of $\delta$.
\item Each set $S_{\delta,i}$ is contained in a strip $R_i$ (i.e. a rectangle of infinite 
length) with width $\delta^{\Cl^{-1} c_0}$.
\item For each $i$ one can find a finite set $F_i \subset \B(0,C)$ of points with
 cardinality 
\be{f-size}
\# F_i \lessapprox \delta^{-C c_0}
\end{equation}
and for each $x \in F_i$ one can associate a collection $A_{i,x,1}, \ldots, 
A_{i,x,M_{i,x}}$ of annuli of thickness $C\delta$, radii which are $\lessapprox 1$ and 
$\delta$-separated, and center $x$ such that
\be{m-size}
M_{i,x} \lessapprox \delta^{-C c_0} \delta^{-1/2}
\end{equation}
for all $x \in F_i$ and
\be{contain}
S_{\delta,i} \subset \bigcup_{x \in F_i} \bigcup_{j=1}^{M_{i,x}} A_{i,x,j}.
\end{equation}
\item We have the estimate 
\be{y-size}
\mu^2(Y_\delta) \lessapprox \delta^{C \Cl^{-1} c_0},
\end{equation}
where
$$ Y_\delta := X_\delta \backslash \bigcup_{i=1}^{N_\delta} S_{\delta,i}^2.$$
\end{itemize}
\end{lemma}

\begin{proof}
Fix $\delta$.  The space of all strips of width $\delta^{\Cl^{-1} c_0}$ which intersect 
$\B(0,C)$ is a two-dimensional manifold, which we can endow with a smooth metric 
$d(\cdot,\cdot)$.  Let $R_1, \ldots, R_{N_\delta}$ be a maximal $C^{-1} 
\delta^{\Cl^{-1} c_0}$-separated subset of this space of strips; note that
$N_\delta \lessapprox \delta^{-2 \Cl^{-1} c_0}$.

For each $x \in K$, let $i(x)$ denote the index $1 \leq i \leq N_\delta$ which maximizes 
the quantity
$$ \mu(\{y \in K: x, y \in R_i \}),$$
and for each $1 \leq i \leq N_\delta$, define $T_{\delta,i}$ to be the set 
$$ T_{\delta,i} := \{ x \in K: d(R_i, R_{i(x)}) \leq C \delta^{\Cl^{-1} c_0} \}.$$
Clearly the sets $T_{\delta,i}$ are contained in $R_i$ and form a finitely overlapping 
cover of $K$.  We shall show the preliminary estimate
\be{t-size}
\mu^2(X_\delta \backslash \bigcup_{i=1}^{N_\delta} T_{\delta,i}^2)
\lessapprox \delta^{C \Cl^{-1} c_0}.
\end{equation}
It suffices to show
$$
\mu^2( \{ (x_0,x_1) \in X_\delta: d(R_{i(x_0)},R_{i(x_1)}) \geq C \delta \} )
\lessapprox \delta^{C \Cl^{-1} c_0}.
$$

From the bounds on $N_\delta$ it suffices to show that
$$
\mu^2( \{ (x_0,x_1) \in X_\delta: i(x_0) = i, i(x_1) = j \} )
\lessapprox \delta^{C \Cl^{-1} c_0}
$$
for all $1 \leq i,j \leq N_\delta$ such that $d(R_i, R_j) \geq C \delta$.

Fix $i$, $j$, and rewrite the above as
$$
\int_{i(x_0) = i} \mu( \{ x_1: (x_0,x_1) \in X_\delta, i(x_1) = j \} )\ d\mu(x_0)
\lessapprox \delta^{C \Cl^{-1} c_0}.
$$
By Cauchy-Schwarz it suffices to show that
$$
\int_{i(x_0) = i} \mu( \{ x_1: (x_0,x_1) \in X_\delta, i(x_1) = j \} )^2\ d\mu(x_0)
\lessapprox \delta^{C \Cl^{-1} c_0}.
$$
By definition of $i(x_0)$, we have
$$ \mu( \{ x_1: (x_0,x_1) \in X_\delta, i(x_1) = j \} )
\leq \mu( \{ x_2: (x_0,x_2) \in X_\delta, i(x_2) = i \} ),$$
so it suffices to show that
$$
\int_{i(x_0) = i} \mu^2( \{ (x_1,x_2): (x_0,x_1), (x_0, x_2) \in X_\delta, i(x_1) = j, 
i(x_2) = i \} )\ d\mu(x_0).
\lessapprox \delta^{C \Cl^{-1} c_0}.
$$
This will obtain if we can show
$$ \mu^3( \{ (x_0,x_1,x_2): x_0,x_2 \in R_i; x_1 \in R_j; (x_0,x_1), (x_0,x_1) \in 
X_\delta \}) \lessapprox \delta^{C \Cl^{-1} c_0}.$$
We first consider the contribution of the case when $|x_1 - x_0| \leq 
\delta^{\Cl^{-1} c_0}$.  In this case we estimate $x_1$ integral by \eqref{frosty} and 
then integrate in the $x_0$ and $x_2$ variables to show that the contribution of this case 
is acceptable.  Similarly we can handle the case $|x_2 - x_0| \leq \delta^{\Cl^{-1} c_0}$. 
 Thus it remains to show that
\begin{align*} \mu^3( \{ (x_0,x_1,x_2): &x_0,x_2 \in R_i; x_1 \in R_j; (x_0,x_1), (x_0,x_1) \in
 X_\delta; |x_1 - x_0|, |x_2 - x_0| > \delta^{\Cl^{-1} c_0} \}) \\
&\lessapprox 
\delta^{C \Cl^{-1} c_0}.
\end{align*}
Suppose $(x_0,x_1,x_2)$ is in the above set.  Since $d(R_i,R_j) > C\delta$, we see from 
elementary geometry that 
$$ |(x_1 - x_0) \wedge (x_2 - x_0)| \geq \delta^{C \Cl^{-1} c_0}.$$
Thus the desired claim follows from \eqref{bilinear-est}, if $\Cl$ is sufficiently large.

The $T_{\delta,i}$ have most of the properties that we desire for $S_{\delta,i}$, but need 
not be covered by a small number of annuli.  To remedy this we shall refine $T_{\delta,i}$ 
slightly.

Fix $j$ and perform the following algorithm.  Initialize $S_{\delta,i}$ to be the 
empty set.  If one has
\be{ash}
\mu( \{ y \in T_{\delta,i} \backslash S_{\delta,i}: (x,y) \in X_\delta \} )
\lessapprox \delta^{Cc_0}
\end{equation}
for all $x \in T_{\delta,i}$, we terminate the algorithm.  Otherwise, we choose an 
$x \in T_{\delta,i}$ for which \eqref{ash} fails, and add the set in \eqref{ash} to 
$S_{\delta,i}$.  We then repeat this algorithm, continuing to enlarge $S_{\delta,i}$ 
until \eqref{ash} is finally satisfied for all $x \in T_{\delta,i}$.

Since each iteration of this algorithm adds a set of measure $\gtrapprox \delta^{Cc_0}$ 
to $S_{\delta,i}$, this algorithm must terminate after at most $\lessapprox \delta^{-Cc_0}$
 steps.  Since $X_\delta$ is a $(\delta, \frac{1}{2} + c_0)_1$ set, we see that each set of 
the form \eqref{ash} can be covered by
annuli $A_{i,x,1}, \ldots, A_{i,x,M_x}$ with width $\delta$, radii $\lessapprox 1$ and 
$C^{-1} \delta$-separated, and center at $x$.  Thus we have the desired covering 
\eqref{covering}.

It remains to prove \eqref{y-size}.  From \eqref{t-size} and the bounds on $N_\delta$ it 
suffices to show that
$$
\mu^2(X_\delta \cap (T_{\delta,i}^2 \backslash S_{\delta,i}^2))
\lessapprox \delta^{C \Cl^{-1} c_0}$$
for all $1 \leq i \leq N_\delta$.  Since
$$ T_{\delta,i}^2 \backslash S_{\delta,i}^2 = T_{\delta,i} \times (T_{\delta,i} \backslash 
S_{\delta,i}) \cup 
 (T_{\delta,i} \backslash S_{\delta,i}) \times T_{\delta,i}$$
it suffices by symmetry to show that
$$
\mu^2(X_\delta \cap (T_{\delta,i} \times (T_{\delta,i} \backslash S_{\delta,i})))
\lessapprox \delta^{C \Cl^{-1} c_0}$$
for all $i$.  But this follows by integrating \eqref{ash} over all $x \in T_{\delta,i}$.
\end{proof}

Henceforth $\Cl$ will be assumed large enough so that the above lemma holds.

From \eqref{y-size} we see that $\sum_\delta \mu^2(Y_\delta) < \infty$, if $\eps$ is chosen 
sufficiently small depending on $c_0$.  From the Borel-Cantelli lemma we thus see that for 
almost every $x$, $y$, the pair $(x,y)$ is contained in only finitely many $Y_\delta$.  
(Here and in the sequel, ``almost every'' is with respect to $\mu$).  Since the $X_\delta$ 
strongly cover $E \times E$, we thus see that $(x,y)$ is contained in infinitely many sets 
of the form $S_{\delta,i}^2$ for hyper-dyadic $\delta$ and $1 \leq i \leq N_\delta$ for 
almost every $x$, $y$.

Suppose $x \in E$, and $\delta$ is a hyper-dyadic number.  Define $d(x,\delta)$ to be the 
smallest hyper-dyadic number $\delta_1$ such that 
\be{dpd}
\delta_1 > \delta^{\Cl/c_0}
\end{equation} 
and such that $x \in S_{\delta_1,i_1}$ for some $1 \leq i_1 \leq N_{\delta_1}$, or 
$d(x,\delta) = +\infty$ if no such $\delta_1$ exists.  From the previous observation we 
thus see that for almost every $x$, $y$, there are infinitely many hyper-dyadic $\delta$, 
$\delta_1$, $\delta_2$ and $i$ such that $x, y \in S_{\delta,i}$, $\delta_1 = d(x,\delta)$, 
and $\delta_2 = d(y,\delta)$.
In particular, we have
$$ \sum_{\delta_1} \sum_{\delta_2} \mu^2 \{ (x,y) \in E \times E: x,y \in S_{\delta,i}, 
\delta_1 = d(x,\delta), \delta_2 = d(y,\delta) \hbox{ for some } \delta, i \} = \infty.$$

The desired contradiction then follows immediately (if $c_0$ is sufficiently small) from 

\begin{lemma}\label{decay}   
For all hyper-dyadic $\delta_1, \delta_2$ we have (if $\Cl$ is chosen appropriately)
\begin{align*}
 \mu^2 \{ (x,y) \in K \times K: &x,y \in S_{\delta,i}, \delta_1 = d(x,\delta), \delta_2 = 
d(y,\delta) \hbox{ for some } \delta, i \} \\
&\leq C_{c_0,\eps}
 \min(\delta_1,\delta_2)^{1/4 - Cc_0}.
\end{align*}
\end{lemma}

The $1/4$ exponent is not optimal, but that is irrelevant for our purposes, since we only 
need the right-hand side to be summable in $\delta_1$, $\delta_2$.

\begin{proof}
Fix $\delta_1, \delta_2$; by symmetry we may assume that $\delta_1 \leq \delta_2$.  
By Fubini's theorem it suffices to show that
\be{concentrate}
\mu \{ x \in K: x,y \in S_{\delta,i}, \delta_1 = d(x,\delta), \delta_2 = d(y,\delta) 
\hbox{ for some } \delta, i \} \leq C_{c_0} \delta_1^{1/4 - Cc_0}
\end{equation}
for all $y \in K$.

Fix $y$.  Since $y$ and $\delta_2$ are fixed, there are significant constraints on the 
number of $\delta$ which can contribute to \eqref{concentrate}.  Indeed, if there are two 
values of $\delta$, say $\delta'$ and $\delta''$, which contribute to \eqref{concentrate}, 
then $\delta'$ cannot exceed ${\delta''}^{\Cl/c_0}$ and $\delta''$ cannot exceed 
${\delta'}^{\Cl/c_0}$, due to the presence of \eqref{dpd} in the definition of 
$d(y,\delta')$, $d(y,\delta'')$.  Because $\delta$ is constrained to be hyper-dyadic, we 
thus see that there are at most $C_{\Cl,c_0,\eps}$ values of $\delta$ which contribute to 
\eqref{concentrate}.  Thus it suffices to show \eqref{concentrate} for a single value of 
$\delta$.  Since the $S_{\delta,i}$ are finitely overlapping as $i$ varies, we see that for 
each $\delta$ there are at most $C$ values of $i$ which contribute to \eqref{concentrate}.  
Hence it suffices to show that
\be{concentrate-2}
\mu \{ x \in K: x \in S_{\delta,i}, \delta_1 = d(x,\delta) \} \leq C_{c_0} 
\delta_1^{1/4 - Cc_0}
\end{equation}
for all $\delta, i$.

Fix $\delta, i$.  We may of course assume that \eqref{dpd} holds, else 
\eqref{concentrate-2}
 is vacuously true.  By definition of $d(x,\delta)$ and the fact that $N_{\delta_1} 
\lessapprox \delta_1^{-C c_0}$ it suffices to show that
$$
\mu \{ x \in K: x \in S_{\delta,i}, x \in S_{\delta_1,i_1} \} \leq C_{c_0} 
\delta_1^{1/4 - Cc_0}
$$
for all $i_1$.

Fix $i_1$.  By \eqref{contain} and \eqref{f-size} it suffices to show that
$$
\sum_{j=1}^{M_{i_1,x_0}} \mu(S_{\delta,i} \cap A_{i_1,x_0,j}) \leq C_{c_0} 
\delta_1^{1/4 - Cc_0}
$$
for all $x_0 \in F$.

Fix $x_0$.  The set $S_{\delta,i}$ is contained in a rectangle of width 
$\delta^{\Cl^{-1} c_0}$, hence contained in a rectangle $R$ of width 
$\delta_1$ by \eqref{dpd}.  Let $r$ denote the distance from $R$ to $x_0$.  
From elementary geometry we see that $R \cap A_{i_1,x_0,j}$ is the union of two sets, 
each of which having diameter at most
$$ \frac{C \delta}{(\delta + |r_j - r|)^{1/2}} $$
where $r_j$ is the outer radius of $A_{i_1,x_0,j}$.  From \eqref{frosty} we thus have
$$
\mu(S_{\delta,i} \cap A_{x_0,j}) \leq C \delta^{1-\eps} (\delta + |r_j - r|)^{1/2-\eps}
$$
If we arrange the $r_j$ in order of distance from $r$, we have $|r_j - r| \geq C j\delta$ 
since the $r_j$ are $\delta$-separated.  Since $\eps \ll c_0$, the claim then follows 
from \eqref{m-size}.
\end{proof}

\end{document}